\documentclass[11pt]{amsart}
\usepackage{amsfonts,amssymb,amsthm}
\usepackage[numbers,sort&compress]{natbib}
\usepackage[leqno]{amsmath}
\usepackage[mathscr]{euscript}
\usepackage{graphicx} 
\usepackage{subfigure}
\usepackage{booktabs}
\usepackage{color}
\usepackage{cases}
\usepackage{mathrsfs}
\usepackage{paralist}
\usepackage{hyperref}
\usepackage{framed}
\usepackage{epstopdf}
\usepackage[ruled, vlined]{algorithm2e}
\usepackage{caption}
\newtheorem {theorem}{Theorem}[section]

\newtheorem {corollary}{Corollary}[section]

\newtheorem {lemma}{Lemma}[section]
\newtheorem {example}{Example}[section]
\newtheorem {definition}{Definition}[section]
\newtheorem {remark}{Remark}[section]

\newcommand{\RR}{{\mathcal R}}
\newcommand{\PP}{{\mathcal P}}
\newcommand{\T}{{\mathcal T}}
\newcommand{\Q}{{\mathcal Q}}
\newcommand{\R}{{\mathbb R}}
\newcommand{\N}{{\mathbb N}}
\newcommand{\Y}{{\mathcal Y}}
\newcommand{\X}{{\mathcal X}}
\newcommand{\A}{{\mathcal A}}
\newcommand{\U}{{\mathcal U}}
\newcommand{\V}{{\mathcal V}}

\newcommand{\cl}{{\rm cl}}
\newcommand{\rank}{{\rm rank}}
\newcommand{\ri}{{\rm ri\;}}

\newcommand{\argmin}{\operatornamewithlimits{argmin}}

\makeatletter
\newcommand{\leqnomode}{\tagsleft@true}
\newcommand{\reqnomode}{\tagsleft@false}
\makeatother

\def\ees{{\accent"5E e}\kern-.385em\raise.2ex\hbox{\char'23}\kern-.08em}
\def\EES{{\accent"5E E}\kern-.5em\raise.8ex\hbox{\char'23 }}
\def\ow{o\kern-.42em\raise.82ex\hbox{
   \vrule width .12em height .0ex depth .075ex \kern-0.16em \char'56}\kern-.07em}
\def\OW{O\kern-.460em\raise1.36ex\hbox{
\vrule width .13em height .0ex depth .075ex \kern-0.16em \char'56}\kern-.07em}

\title[Utopia point method based robust vector polynomial optimization scheme]
{A utopia point method-based robust vector polynomial optimization scheme$^*$}

\thanks{$^*$The authors (in alphabetic order) contributed equally to this work}

\author[Tianyi Han]{Tianyi Han}
\address[Tianyi Han]{School of Mathematical Sciences, Shanghai Jiao Tong University, Shanghai 200240, China}
\email{18846073679@163.com; hantianyi@sjtu.edu.cn}

\author{Liguo Jiao}
\address[Liguo Jiao]{Academy for Advanced Interdisciplinary Studies, Northeast Normal University, Changchun 130024, Jilin Province, China.}
\email{hanchezi@163.com; jiaolg356@nenu.edu.cn}

\author{Jae Hyoung Lee}
\address[Jae Hyoung Lee]{Department of Applied Mathematics, Pukyong National University, Busan 48513, Korea.}
\email{mc7558@naver.com}

\author{Junping Yin$^{**}$}
\address[Junping Yin]{Institute of Applied Physics and Computational Mathematics,Beijing 100094, China; Shanghai Zhangjiang Academy of Mathematics, Shanghai 201203, China; Academy for Advanced Interdisciplinary Studies, Northeast Normal University, Changchun 130024, Jilin Province, China.}
\email{yinjp829829@126.com}

\thanks{$^{**}$Corresponding Author}



\date{\today}

\begin{document}

\maketitle

\begin{abstract}
In this paper, we focus on a class of robust vector polynomial optimization problems (RVPOP in short) without any convex assumptions. 
By combining/improving the utopia point method (a nonlinear scalarization) for vector optimization and ``joint+marginal" relaxation method for polynomial optimization, we solve the RVPOP successfully. 
Both theoratical and computational aspects are considered.
\end{abstract}


\section{Introduction}

By vector optimization, we mean minimizing simultaneously several objective functions over a constraint set as usual.
In contrast with numerical optimization, it is well-known that one may not be possible to find a single point that simultaneously minimizes all the objective functions in vector optimization.
In this case, we always look for some ``best preferred" solution(s), by which, we mean a solution that cannot be improved in one objective function without deteriorating their performance in at least one of the rest.
This concept is nowadays known as the Pareto/efficient solution in vector optimization; see, e.g., \cite{Chankong1983,Luc1989,Ehrgott2005,Jahn2011,Sawaragi1985}.

In most real-world problems, the related parameters and coefficients usually depend on forecasting, approximations, or may vary due to disturbances or other uncertainties, such problems with data uncertainty are studied as uncertain optimization problems.
As a rather powerful deterministic approach for studying uncertain optimization problems, the robust optimization can associate such problems with its robust counterpart; see, e.g., \cite{Ben-Tal2009,Ben-Tal2002,Ben-Tal2008,Bertsimas2011}.

We shall mention here that, the robust optimization usually deals with uncertain data, while vector optimization studies minimizing simultaneously several objective functions, the problem combining ``robust" with ``vector" optimization is known as {\it robust vector optimization problem} (RVOP in short).
Indeed, the most important issue in RVOPs is to carefully define what a robust solution to an uncertain vector problem is.
This issue is extremely difficult, since not only the notion of {\it robustness} depends on the chosen concept but also {\it efficiency} shall be defined with respect to certain ordering cones;
in order to solve this difficulty, different concepts of robust efficiency have been given and studied in the literature; see, e.g., \cite{Bitran1980,Soyster1973,Beck2009,kuroiwa2012robust}.
Among them, we are concerning with the one (called minimax robust efficiency) proposed by Kuroiwa and Lee \cite{kuroiwa2012robust} from a set-valued point of view, in which case, the uncertainty is assumed to be independent for each objective function.

\subsection{Problem formulations and solution concepts}
In this paper, we are interested in a RVOP with {\it polynomial data}.
More precisely, among other things, we aim to solve it, i.e., finding its properly robust Pareto solutions (see Definition~\ref{properParetosolu}).
Mathematically, let us consider the following class of uncertain vector polynomial optimization problems:
\begin{align*}\label{UVOP}
\min~&f(x,u) \tag{UVPOP}\\
{\rm s.t.}~~~&g_j(x,v) \geq0, ~ j=1,\ldots,m,
\end{align*}
where $f(x,u):= \left(f_1(x,u),\ldots,f_l(x,u)\right),$ $f_i\colon \R^n\times\R^p \to \R,$ $i=1,\ldots,l,$ and $g_j\colon \R^n\times\R^q \to \R,$ $j=1,\ldots,m,$ are real polynomials. 
The uncertainty set consisting of $u$ is assumed to be a compact basic semi-algebraic set described by
\begin{equation*}
\U:=\{u=(u_1,\ldots,u_p)\in\R^p\colon \hat h_k(u) \ge 0,~k=1,\ldots,s_1\},
\end{equation*}
where $\hat h_k:\R^p\to\R,~k=1,\ldots,s_1$ are polynomials; and the uncertainty set consisting of $v$ is also a compact basic semi-algebraic set given by
\begin{equation*}
	\V:=\{v=(v_1,\ldots,v_q)\in\R^q\colon \bar h_k(v) \ge 0,~k=1,\ldots,s_2\},
\end{equation*}
where $\bar h_k:\R^q \to \R,~k=1,\ldots,s_2$ are also polynomials.

In this work, we mainly consider a pessimistic situation;
roughly speaking, we aim to find feasible solutions that are less sensitive to the perturbations of uncertainty parameters $u$ and $v.$
Therefore, following the definitions proposed by Kuroiwa and Lee~\cite{kuroiwa2012robust}, we consider the worst cases of both objective functions and constraint sets, that is, the robust counterpart to the uncertain vector polynomial optimization problem~\eqref{UVOP},
\begin{align*}\label{RVOP}
\min_{x}~~&\sup_{u\in \U} f(x,u)\tag{RVPOP}\\
{\rm s.t.}~~&\inf_{v\in \V}g_j(x,v)\ge0,~j=1,\ldots,m.
\end{align*}

Below, we recall the concept of robust Pareto solution to the problem~\eqref{UVOP}, note first that a point $x\in\R^n$ is said to be \textit{robust feasible} to the problem~\eqref{UVOP}, if $\inf\limits_{v\in \V}g_j(x,v)\ge0,~j=1,\ldots,m,$ and all the robust feasible solution forms the {\it robust feasible solution set} denoted by
\begin{align}\label{ro_fe_so_se}
\Omega:=\left\{x \in \R^n \colon \inf\limits_{v\in \V}g_j(x,v)\ge0,\ j=1,\ldots,m\right\}.
\end{align}
Throughout this paper, we make the following blanket assumption.
\begin{framed}
	\begin{itemize}
		\item[{\bf (H1)}] We always assume the robust feasible solution set $\Omega$ is nonempty and compact.
	\end{itemize}
\end{framed}

\begin{remark}\label{H2}
{\rm
By assumption {\bf (H1)}, it is reasonable for us to assume that there exisits a simple compact semi-algebraic set (e.g., a box or an ellipsoid)
\begin{align}\label{simpleset}
X:=\{x\in\R^n\colon \tilde h_k(x)\ge0,k=1,\ldots,s_3\}
\end{align}
containing all feasible $x$, where $\tilde h_k:\R^n\to\R,~k=1,\ldots,s_3$ are polynomials. Namely, it holds that
\begin{equation*}
	\{x\in\R^n\colon \inf_{v\in \V} g_j(x,v) \geq0, ~ j=1,\ldots,m\}\subseteq X.
\end{equation*}
}
\end{remark}

\begin{definition}[robust weak Pareto solution]{\rm \cite{kuroiwa2012robust}
	A point $\hat x \in \Omega$ is said to be a \emph{robust weak Pareto solution} to the problem~\eqref{UVOP}, if $\hat{x}$ is robust feasible and
	\begin{equation*}
	\sup_{u\in \U}f(x,u)-\sup_{u\in \U}f(\hat x,u)\not\in-\R^l_> \ \  \textrm{for all}\ x \in \Omega,
	\end{equation*}
where $\R^l_>:=\left\{y:=(y_1,\ldots, y_l) \colon y_i >0,\ i = 1, \ldots, l\right\}$ stands for the positive orthant of $\R^l.$
}\end{definition}

\begin{definition}[robust Pareto solution]\label{Paretosolu}
{\rm \cite{kuroiwa2012robust}	A point $\hat x\in \Omega$ is said to be a \emph{robust Pareto solution} to the problem~\eqref{UVOP}, if $\hat{x}$ is robust feasible and
	\begin{equation*}
		\sup_{u\in \U}f(x,u)-\sup_{u\in \U}f(\hat x,u)\not\in-\R^l_{\geq}\setminus\{0\} \ \ \textrm{for all}\ x \in \Omega.
	\end{equation*}
Moreover, if $\hat x$ is a robust Pareto solution to the problem~\eqref{UVOP}, then $f(\hat x) \in \R^l$ is called a {\it robust Pareto value} to the problem~\eqref{UVOP}, and all robust Pareto values form the {\it robust Pareto curve} (or robust Pareto front).
}\end{definition}

\begin{definition}[properly robust Pareto solution]\label{properParetosolu}
{\rm \cite{kuroiwa2012robust}	
A point $\hat x\in \Omega$ is said to be a \emph{properly robust Pareto solution} to the problem~\eqref{UVOP} if $\hat{x}$ is robust feasible and there exsits $M>0$ such that for all indices $i$ and robust feasible points $x$ satisfying $\sup\limits_{u\in \U}f_i(x,u)<\sup\limits_{u\in \U}f_i(\hat x,u)$, there exists an index $i_0$ with $\sup\limits_{u\in \U}f_{i_0}(\hat x,u)<\sup\limits_{u\in \U}f_{i_0}(x,u)$ such that
	\begin{equation*}
		\frac{\sup_{u\in \U}f_i(\hat x,u)-\sup_{u\in \U}f_i(x,u)}{\sup_{u\in \U}f_{i_0}(x,u)-\sup_{u\in \U}f_{i_0}(\hat x,u)}\le M.
	\end{equation*}
}\end{definition}

\subsection{The utopia point method scheme}
It is known that a vector optimization problem can be solved by solving related single objective optimization problems.
Such methods are called scalarization approaches, and there are many types of scalarization approaches (see \cite{Chankong1983,Luc1989,Ehrgott2005,Lee2021JOGO}); for instance, the linear scalarization methods including the weighted sum method and the hybrid method.
It is worth noting that, the above mentioned two methods have been proved to be effective for finding (robust) Pareto solutions to vector optimization problem with {\it convex polynomials}; see, e.g., \cite{Chuong2018,Chuong2022,Jiao2019,Jiao2020,Lee2018,Lee2021,Lee2021JOGO}.
However, for solving the problem~\ref{RVOP}, a {\it nonconvex} case, such linear scalarization methods are no longer valid.
This deeply motivates us to find/integrate other (nonlinear scalarization) methods for solving the problem~\ref{RVOP}.

In this part, we shortly depict the so-called utopia point method --- a nonlinear scalarization approach---for solving the problem \eqref{RVOP}.
Firstly, inspired by \cite[Definitions 2.22 \& 4.23]{Ehrgott2005}, we give the following concepts for the problem \eqref{RVOP}.
\begin{definition}{\rm
		The point $y^I=(y^I_1, \ldots, y^I_l)$ given by
		\begin{equation}\label{ideal}
		y^I_i:=\min_{x\in \Omega}\sup_{u\in \U}f_i(x,u), \ ~i=1,\ldots,l,
		\end{equation}
		where $\Omega$ is defined as in \eqref{ro_fe_so_se},	is called the {\it ideal point} of the problem \eqref{RVOP}.
}\end{definition}

\begin{definition}\label{utopia}{\rm
		The point $y^U:=y^I-\epsilon$ with $\epsilon\in\R^l_>$ has positive components, is called a {\it utopia point} of the problem \eqref{RVOP}.
}\end{definition}

The figure below demonstrates the difference between the ideal point and a utopia point for general vector optimization problems in the objective space .

\begin{figure}[htbp]
	\includegraphics*[width=2in]{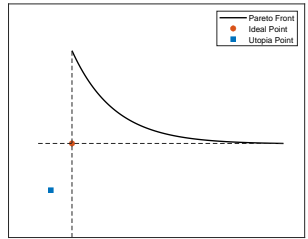}
	\caption{The ideal point and a utopia point}\label{explain}
\end{figure}

Denoted by $\lambda:=(\lambda_1,\ldots,\lambda_l)$ certain weighting parameter in $\R^l$ and $\|\cdot\|_p$ the $\ell_p$-norm defined in $\R^l$, let $\ri\Lambda:=\left\{\lambda\in\R^l_>\colon \sum_{i=1}^{l}\lambda_i=1\right\}$ stand for the relative interior of the unit simplex.
In this paper, we are going to deploy an approximation method (i.e., utopia point method) for solving the problem \eqref{RVOP}, which leads to the following optimization problem,
\begin{equation}\label{approx}
	\min_{x\in \Omega}\|\lambda\odot(\sup_{u\in \U}f(x,u)-y^U)\|^p_p, \tag{${\rm P^{utopia}}$}
\end{equation}
where $a \odot b :=(a_1b_1,\ldots,a_lb_l)$ refers to the Hadamard product for $a, b \in \R^l$ and $\lambda\in\ri\Lambda$.
Indeed, by solving the problem~\eqref{approx}, we can extract the (properly) robust Pareto solutions to the problem~\eqref{UVOP}, this is also the main findings in this work; see Section~\ref{sect::4} in detail.

The reason why we employ the utopia point in \eqref{approx} rather than the ideal point basically owes to the following three aspects:
\begin{itemize}
	\item [{\rm A.}] Usually, it is not easy to solve \eqref{ideal} exactly to obtain an ideal point in the context that $f_i(x,u)$ is nonconvex in $x$ and nonconcave in $u$, whereas estimating a lower bound of the ideal point with a special relaxation method is much easier; one may compare with the results in the work~\cite{GuoJiao2022arXiv}.
	\item [{\rm B.}] Particularly, if we drop the uncertainty parameter in \eqref{UVOP}, it will be converted into a vector polynomial optimization problem reading as
	\begin{align*}
		\min\ &f(x):=\left(f_1(x), \ldots, f_l(x)\right)\\
		{\rm s.t.}\ &g_j(x)\ge0,~j=1,\ldots,m,
	\end{align*}
	thus we reduce the problem~\eqref{approx} into
	\begin{equation}\label{approxp}
		\min \left\{\|\lambda\odot(f(x)-y^U)\|^p_p\colon g_j(x)\ge0,~j=1,\ldots,m \right\}.
	\end{equation}
It has been proven that the closure of Pareto front can be characterized by \eqref{approxp}, whereas it doesn't hold if we replace the utopia point $y^U$ in \eqref{approxp} with the ideal point $y^I$.
More details can be found in \cite[Chapter 4.5]{Ehrgott2005}.

\item [{\rm C.}] Choosing a utopia point covers many special cases, whereas the ideal point method fails to cover such cases.
For instance, when the $\ell_p$-norm in \eqref{approx} is set to be the $\ell_\infty$-norm, it coincides with the well-known Chebyshev scalarization method,
\begin{equation}\label{chebyshev}
	\min_{x\in \Omega}\max_{1\le i\le l}\ \lambda_i(\sup_{u\in \U}f_i(x,u)-y_i^U).
\end{equation}
Furthermore, we have a robust feasible solution $\hat{x}$ is a robust weak Pareto solution if and only if $\hat{x}$ is an optimal solution to the problem~\eqref{chebyshev}; see, e.g., \cite[Theorem 4.24]{Ehrgott2005}.
Besides, compared with the Chebyshev scalarization method, the utopia point method shares the benefit of computing the more valuable properly robust Pareto solutions and values, however it has to bear the drawback of higher computational complexity when inceasing the order of $\ell_p$-norm.
To relieve this problem, sparsity should be taken into consideration, but this is out of the scope of this paper.
\end{itemize}

\subsection{Contributions}
In the present paper, we mainly make the following {\it contributions} to the robust vector polynomial optimization problem~\eqref{RVOP}.

\begin{itemize}
\item [{\rm (i)}] Comparing with our previous works \cite{Jiao2019,Jiao2020,Lee2018,Lee2021,Lee2021JOGO}, we don't the considered polynomials to be convex, hence the usual linear scalarization methods are not valid.
\item [{\rm (ii)}] Besides, even though $f_i(x,u)$ and $g_j(x,v)$ $i=1,\ldots,l,~j=1,\ldots,m$ are polynomials, the supremum and infimum of which are no longer polynomials, consequently we shall deploy new ``joint+marginal"-based relaxation methods to design the algorithms.
\item [{\rm (iii)}] We develop the utopia point method --- a novel nonlinear scalarization approach---for solving the problem \eqref{RVOP}.
\end{itemize}

The rest of the paper is organized as follows.
Section~\ref{sect::2} gives some basic notations and preliminaries that will be used in this paper.
In section \ref{sect::3} we deploy ``joint+marginal" relaxation method and obtain the upper/lower bound polynomials for the supreme/infimum of the objective/constraint functions, and demonstrate the convergence results.
In section \ref{sect::4} we first establish the concepts of various approximation sets, then we give the convergence results for them.
Section \ref{sect::5} gives the procedure of our algorithm and some illustrative paradigms.
Section \ref{sect::6} summarizes the obtained results.

\section{Preliminaries}\label{sect::2}

We begin this section by fixing some notations and preliminaries; see \cite{Lasserre2010,Lasserre2015} for more details.
We suppose $1 \leq n \in \N$ ($\N$ is the set of non-negative integers) and abbreviate $(x_1, x_2, \ldots, x_n)$ by $x.$
The Euclidean space $\R^n$ is equipped with the usual Euclidean norm $\| \cdot \|$ and the inner product $\langle \cdot, \cdot \rangle,$ here $\langle x, y \rangle := x^Ty = \sum_{i = 1}^n x_iy_i,$ for any $x, y \in \R^n.$
The non-negative orthant of $\R^n$ is denoted by $\R_{+}^n.$

Let $\R[x]$ stand for the space of all real polynomials in the variable $x$, and $\R[x]_d$ be the space of all real polynomials in the variable $x$ with degree at most $d$. The degree of a polynomial $f$ is denoted by $\deg f.$
We say that a real polynomial $f$ is sum-of-squares (SOS) if there exist real polynomials $q_l,$ $l = 1,\ldots,r,$ such that $f =\sum_{l=1}^{r}q_{l}^2.$
The set consisting of all sum-of-squares real polynomials is denoted by $\Sigma[x],$ and the set consisting of all sum-of-squares real polynomials with degree at most $d$ is denoted by $\Sigma[x]_{2d}$, which is a closed cone.

Given $m$ polynomials $\{g_1,\ldots,g_m\} \subset \R[x],$ a {\it quadratic module} generated by the tuple $g:=(g_1,\ldots,g_m)$ is defined by
\begin{equation*}
\mathcal{Q}(g) := \left\{\sigma_0 + \sum_{i=1}^m\sigma_ig_i: \sigma_i \in \Sigma[x], \ i = 0, 1, \ldots, m\right\}.
\end{equation*}
Define $g_0\equiv1$, the $2d$-th truncation of $\Q(g)$ is
\begin{equation*}
	\Q_{2d}(g) := \left\{\sigma_0 + \sum_{i=1}^m\sigma_ig_i\colon \sigma_i \in \Sigma_{2d-\deg g_i}[x], \ i = 0, 1, \ldots, m\right\}.
\end{equation*}

Let $\PP[K]$ stand for the space of all non-negative polynomials on the set $K,$ and denoted by $\PP_d[K]:=\PP[K]\bigcap\R[x]_d$ the space of all non-negative polynomials on the set $K$ with degree at most $d.$ 
When $K=\{x\in\R^n: g_i(x) \geq 0, \ i=1,\ldots,m\}$, it holds for all degree $d$ that
\begin{equation*}
	\Q_d(g)\subseteq\PP_d[K].
\end{equation*}

\begin{definition}{\rm
We say the set $\mathcal{Q}(g)$ satisfies the \emph{Archimedean property} if there exists $p\in \mathcal{Q}(g)$ such that the set $\{x \in\R^n: p(x) \geq 0\}$ is compact.
}\end{definition}

It can be justified that $K$ must be compact if $\Q(g)$ is Archimedean, but the converse is not necessarily true. The following result describes the celebrated representation of positive polynomials over a semi-algebraic set when the quadratic module $\mathcal{Q}(g)$ is Archimedean.
\begin{lemma}[Putinar's Positivstellensatz]
{\rm \cite{Putinar1993}
Let $f \in \R[x]$ and $g_i \in \R[x],$ $i=1,\ldots,m.$
Suppose that $\mathcal{Q}(g)$ is Archimedean.
If $f$ is strictly positive on $K:=\{x\in\R^n: g_i(x) \geq 0, \ i=1,\ldots,m\}\neq\emptyset,$ then $f \in\mathcal{Q}(g).$
}\end{lemma}

For a multi-index $\alpha:=(\alpha_1,\ldots,\alpha_n)\in\N^n,$ let us denote $|\alpha|:=\sum_{i=1}^n\alpha_i$ and $\N^n_d:= \{\alpha \in \N^n :  |\alpha|\leq d\}.$
The notation $x^\alpha$ denotes the monomial $x_1^{\alpha_1}\cdots x_n^{\alpha_n}.$
The canonical basis of $\R[x]_d$ is denoted by
\begin{equation*}
v_d(x):=(x^\alpha)_{\alpha\in\N^n_{d}}=(1,x_1,\ldots,x_n,x_1^2,x_1x_2,\ldots,x_n^2,\ldots,x_1^d,\ldots,x_n^d)^T,
\end{equation*}
which has dimension $s(d):=\left( \substack{ n+d \\ d }\right).$ It holds that for all $f\in\Sigma_{2d}[x]$ there exists a semi-definite symmetirc matrix $Q\in\R^{s(d)\times s(d)}$ such that $f(x)=v_d(x)^TQv_d(x)$.

Given a truncated sequence of degree no more than $d$
\begin{equation*}
	y=(y_\alpha)_{\alpha\in\N^n_d},
\end{equation*}
the corresponding Riesz functional $L_y:\R_d[x]\to\R$ is given by
\begin{equation*}
	f\left(=\sum_{\alpha\in\N^n_d}f_\alpha x^\alpha\right)\mapsto\sum_{\alpha\in\N^n_d}f_\alpha y_\alpha\left(=:L_y(f)\right).
\end{equation*}
Given an $s(2d)$-sequence $y:= (y_\alpha)_{\alpha\in\N^n_{2d}}$ with $y_0=1$, we define $d$-th order moment matrix as
\begin{equation*}
	M_d(y):=L_y\left(v_d(x)v_d(x)^T\right),
\end{equation*}
for a polynomial $g\in\R[x]$, the $d$-th order localizing matirx is
\begin{equation*}
	M_d(gy):=L_y\left(g(x)v_d(x)v_d(x)^T\right).
\end{equation*}

Recall that $\PP[K]$ denotes the cone of all nonnegtive polynomial in $K,$ we define the dual cone of $\PP[K]$ by Riesz functional
\begin{equation*}
	\PP^*[K]:=\left\{(y_\alpha)_{\alpha\in\N^n}\colon L_y(f)\ge0,~f\in\PP[K]\right\},
\end{equation*}
if $K$ is the semi-algebraic set defined as above and Archimedean condition is satisfied, then it holds that $\PP^*[K]=\Q^*(g)$. While for the truncated case, one only has
\begin{equation*}
	\PP^*_d[K]\subseteq\Q^*_d(g).
\end{equation*}

Assume that the largest degree within the constraint polynomials $g_j$ is no more than $2d,$ define the $2d$-tms cone
\begin{equation*}
	\T_{2d}(g):=\{y\in\N^n_{2d}\colon M_d(g_jy)\succeq0,~j=0,\ldots,m\},
\end{equation*}
where ``tms" is the abbreviation of \textit{truncated moment sequence}. It can be verified that (see \cite{nie2015linear})
\begin{equation*}
	\Q^*_{2d}(g)=\T_{2d}(g).
\end{equation*}
A tms $y$ is said to admit a representing measure $\mu$ supported on $K$ if $y_\alpha=\int_Kx^\alpha d\mu$ for all $\alpha\in\N^n$, such a measure is called an $K$-representing measure for $y.$ 
The set of all tms admitting $K$-representing measure is given by the moment cone
\begin{equation*}
	\RR(K):=\left\{y=(y_\alpha)_{\N^n}\colon \exists \mu\ \textrm{supported in $K$ s.t.}\ y_\alpha=\int_Kx^\alpha d\mu,\ \forall\alpha\in\N^n\right\},
\end{equation*}
for $d$-tms $y=(y_\alpha)_{\N^n_d}$ admitting $K$-representing measure, we define the corresponding moment cone as $\RR_d(K)$.

For all $f\in\PP_d[K]$ and $y\in\RR_d(K)$, there exists a representing measure of $y$ such that
\begin{equation*}
	L_y(f)=\sum_{\alpha\in\N^n_d}f_\alpha y_\alpha=\int_K f(x)d\mu\ge0,
\end{equation*}
this implies that $\RR^*_d[K]=\PP_d(K)$, moreover if $K$ is compact one has $\PP^*_d[K]=\RR_d(K)$, otherwise $\PP^*_d[K]=\cl(\RR_d(K))$ \cite{nie2015linear}.

By now, if $K$ is the semi-algebraic set defined as above and Archimedean condition is satisfied, we have the following relationship:
\begin{equation*}
	\RR_{2d}(g)=\PP^*_{2d}[K]\subseteq\Q^*_{2d}(g)=\T_{2d}(g).
\end{equation*}

\begin{definition}[flat truncation]
{\rm \cite{nie2013certifying}
	Let $K=\{x\in\R^n: g_i(x) \geq 0, \ i=1,\ldots,m\}$, for an integer $d\ge \deg g:=\max\limits_{j=1,\ldots,m}\deg g_j$, a tms $y\in\RR_{2d}(K)$ is called {\it flat truncated} if there exists $t\le d$ such that
	\begin{equation*}
		\rank(M_{t-d_0}(y))=\rank(M_t(y))=:r,
	\end{equation*}
	where $d_0=\lceil \deg g/2 \rceil$.
}\end{definition}

If a tms $y$ is flat, then the corresponding measure $\mu$ is unique and $r$-atomic \cite{nie2013certifying,nie2014optimality}.

\section{The ``joint+marginal" relaxation method}\label{sect::3}

In this section, we are going to tackle the intermediate problem of \eqref{RVOP}. Consider the following parametric polynomial optimization problems
\begin{align}
	F_i(x):=&\sup_{u\in \U}~f_i(x,u),~\forall x\in X,~i=1,\ldots,l, \label{paraf} \\
	G_j(x):=&\inf_{v\in \V}~g_j(x,v),~\forall x\in X,~j=1,\ldots,m, \label{parag}
\end{align}
where $X$ is defined as in \eqref{simpleset}.
Once $F_i(x)$ and $G_j(x)$, the optimal value functions of parametric polynomial optimization problems \eqref{paraf} and \eqref{parag} respectively, are found for $i=1,\ldots,l$ and $j=1,\ldots,m$, we can convert the robust vector polynomial optimization problem \eqref{RVOP} into the following vector optimization problem:
\begin{align*}\label{VP}
	\min_{x}~&F(x):=\left(F_1(x),\ldots,F_l(x)\right)\tag{VP}\\
	{\rm s.t.}~&G_j(x)\ge0,~j=1,\ldots,m.
\end{align*}
Observe that, the supremum and infimum of a family of polynomials are no longer polynomials, and this may be the biggest obstacle when solving \eqref{VP}. 
Nevertheless, inspired by \cite{lasserre2011min}, we can approximate the optimal value function of the parametric polynomial optimization problem \eqref{paraf} with "joint+marginal" method proposed in \cite{lasserre2010joint+}, which generates a sequence of polynomials optimal upper/lower approximations with strong convergence in some measure to the supremum/infimum value function. 
With this approximation methodology we drop the robustness and obtain a sequence of vector polynomial optimization problems, which can be solved with Lasserre's hierarchy. 
It is worth mentioning that we can obtain the approximation as closed as desired theoretically without any discretization scheme.

Let $\varphi$ be a strictly positive Borel probabilty measure on $X,$ and is absolutely continuous with respect to the Lebesgue measure $m$ defined on $\R^n$. 
We employ the measure presented in \cite{lasserre2010joint+}, choose $\varphi$ the probability measure uniformly distributed on $X,$
\begin{equation*}
	\varphi(B):=\left(\int_Xdm\right)^{-1}\int_Bdm,~~\forall B\in\mathcal{B}(X).
\end{equation*}
Therefore we are able to compute all moments with respect to $\varphi$ on $X$, denoted by $\gamma=(\gamma)_\alpha$, $\alpha\in\N^n$, with
\begin{equation*}
	\gamma_\alpha:=\int_Xx^\alpha d\varphi(x),
\end{equation*}
that is the reason why $X$ is chosen to be a simple set.

Now, consider the optimization problems
\begin{align}\label{marginalf}
	\min_{\mu^i}~&\sum_{\alpha\in\N^n}\mu^i_\alpha\gamma_\alpha \\
				 {\rm s.t.}~&\sum_{\alpha\in\N^n}\mu^i_\alpha x^\alpha-f_i(x,u)\in\PP[X\times \U],\notag
\end{align}
and
\begin{align}\label{marginalg}
	\max_{\nu^j}~&\sum_{\alpha\in\N^n}\nu^j_\alpha\gamma_\alpha \\
	{\rm s.t.}~&g_j(x,v)-\sum_{\alpha\in\N^n}\nu^j_\alpha x^\alpha\in\PP[X\times \V].\notag
\end{align}
To proceed, by \cite{lasserre2010joint+}, we have
\begin{equation*}
	\sum_{\alpha\in\N^n}(\mu^i_\alpha)^*\gamma_\alpha=\int_{x\in X}F_i(x)d\varphi(x),
\end{equation*}
and
\begin{equation*}
\sum_{\alpha\in\N^n}(\nu^j_\alpha)^*\gamma_\alpha=\int_{x\in X}G_j(x)d\varphi(x),
\end{equation*}
where $(\mu^i_\alpha)^*$ and $(\nu^j_\alpha)^*$ achieve the minimum of \eqref{marginalf} and the maximum of \eqref{marginalg}, respectively, for $i=1,\ldots,l$ and $j=1,\ldots,m$.

For every $d^i\in\N$, we have semidefinite programming relaxation of \eqref{marginalf} with relaxation order $(d^i)^\prime\ge\max\left\{\lceil (d^i/2 \rceil,\lceil (\deg f_i)/2 \rceil,\lceil \max\limits_{k=1,\ldots,s_1}(\deg \hat h_k)/2 \rceil, \lceil \max\limits_{k=1,\ldots,s_3}(\deg \tilde h_k)/2 \rceil \right\}$ when $X\times \U$ is Archimedean:
\begin{align}\label{sdpmarginalf}
	\min_{\mu^i}~&\sum_{\alpha\in\N^n_{d^i}}\mu^i_\alpha\gamma_\alpha \\
	{\rm s.t.}~&\sum_{\alpha\in\N^n_{d^i}}\mu^i_\alpha x^\alpha-f_i(x,u)\in\Q_{2(d^i)^\prime}(\tilde h,\hat h).\notag
\end{align}
For all feasible $\mu^i$ of \eqref{sdpmarginalf}, let $F^{d^i}(x):=\sum\limits_{\alpha\in\N^n_{d^i}}\mu^i_\alpha x^\alpha\in\R_{2d^i}[x]$, it holds that $F^{d^i}(x)\ge F_i(x)$ in $X$ for $i=1,\ldots,l,$ since
\begin{equation*}
	F^{d^i}(x)-f_i(x,u)\in\Q_{2(d^i)^\prime}(\tilde h,\hat h).
\end{equation*}
On the other hand
\begin{equation*}
	\sum_{\alpha\in\N^n_{d^i}}\mu^i_\alpha\gamma_\alpha=\int_{x\in X}F^{d^i}(x)d\varphi(x),
\end{equation*}
therefore \eqref{sdpmarginalf} generates an upper bound of $F_i(x)$.

For every $e^j\in\N$, similar to \eqref{sdpmarginalf}, we have semidefinite programming relaxation of \eqref{marginalg} with relaxation order $(e^j)^\prime\ge \max\left\{\lceil (e^j/2 \rceil,\lceil (\deg g_j)/2 \rceil,\lceil \max\limits_{k=1,\ldots,s_2}(\deg \bar h_k)/2 \rceil, \lceil \max\limits_{k=1,\ldots,s_3}(\deg \tilde h_k)/2 \rceil \right\}$, when $X\times \V$ is Archimedean:
\begin{align}\label{sdpmarginalg}
\max_{\nu^j}~&\sum_{\alpha\in\N^n_{e^j}}\nu^j_\alpha\gamma_\alpha \\
{\rm s.t.}~&g_j(x,v)-\sum_{\alpha\in\N^n_{e^j}}\nu^j_\alpha x^\alpha\in\Q_{2(e^j)^\prime}(\tilde h,\bar h).\notag
\end{align}
Moreover, denoted by $G^{e^j}(x):=\sum\limits_{\alpha\in\N^n_{e^j}}\nu^j_\alpha x^\alpha\in\R_{e^j}[x]$ the lower bound of $G_j(x)$, we have the following convergence result, which tells us that we can approximate the optimal functions $F_i(x)$ and $G_j(x)$ with polynomial in a strong sense, and even in a stronger sense with piecewise polynomial.

\begin{theorem}\label{convergence}{\rm (c.f. \cite{lasserre2010joint+})}
	Suppose that $\mu^i$ and $\nu^j$ are optimal or nearly optimal solution to \eqref{sdpmarginalf} and \eqref{sdpmarginalg} respectively$,$ let $F^{d^i}(x):=\sum_{\alpha\in\N^n_{d^i}}\mu^i_\alpha x^\alpha\in\R_{d^i}[X]$ and $G^{e^j}(x):=\sum_{\alpha\in\N^n_{e^j}}\nu^j_\alpha x^\alpha\in\R_{e^j}[x],$ it holds that$:$
	\begin{itemize}
		\item [{\rm (a)}]$||F^{d^i}-F_i||_{L_1(X,\varphi)}\overset{d^i\to\infty}{\longrightarrow}0$ and $||G^{e^j}-G_j||_{L_1(X,\varphi)}\overset{e^j\to\infty}{\longrightarrow}0$.
		\item [{\rm (b)}]Let $\tilde{F}^{d^i}(x):=\min\limits_{d\le d^i}\{F^d(x)\}$ and $\tilde{G}^{e^j}(x):=\max\limits_{e\le e^j}\{G^e(x)\},$ then $\tilde{F}^{d^i}(x)\downarrow F_i(x)$ and $\tilde{G}^{e^j}(x)\uparrow G_j(x)$ for all $x\in X,$ besides $\tilde{F}^{d^i}(x)\to F_i(x)$ and $\tilde{G}^{e^j}(x)\to G_j(x)$ $\varphi$-almost uniformly in $X$.
	\end{itemize}
\end{theorem}

Let $z_{(d^i)^\prime}$ be a $2(d^i)^\prime$-th truncated moment sequence defined on $\R^n\times\R^p$, which means
\begin{equation*}
	z^{\alpha\beta}_{(d^i)^\prime}:=x^\alpha u^\beta,~~\forall \alpha,\beta\in\N^n_{2(d^i)^\prime},~|\alpha|+|\beta|\le s(2(d^i)^\prime).
\end{equation*}
Then the dual problem of \eqref{sdpmarginalf} reads as
\begin{align}\label{dualsdpmarginal}
	\max_{z_{(d^i)^\prime}}~&L_{z_{(d^i)^\prime}}(f_i)\\
	{\rm s.t.}~&z_{(d^i)^\prime}\in\T_{2(d^i)^\prime}(\tilde h,\hat h)\notag\\
		~&L_{z_{(d^i)^\prime}}(x^\alpha)=\gamma_\alpha,~\forall\alpha\in\N^n_{d^i}\notag\\
		~&L_{z_{(d^i)^\prime}}(x^\alpha)=0,~s(d^i)\le|\alpha|\le s(2(d^i)^\prime),\notag
\end{align}
which can also be solved by semidefinite programming. 
For all feasible $z_{(d^i)^\prime}$ of \eqref{dualsdpmarginal}, let $u:=(z^{0e_1}_{(d^i)^\prime},\ldots,z^{0e_p}_{(d^i)^\prime})$, then $u\in \U$ since $z_{(d^i)^\prime}\in\T_{2(d^i)^\prime}(\tilde h,\hat h)$. 
On the other hand $L_{z_{(d^i)^\prime}}(x^\alpha)=\gamma_\alpha,~\forall\alpha\in\N^n_{2d^i}$, therefore $L_{z_{(d^i)^\prime}}(f_i)\ge \int_{x\in X}f_i(x,u)d\varphi(x)$. Moreover,

\begin{theorem}\label{momentconvergence}{\rm (c.f. \cite{lasserre2010joint+})}
	Suppose z is an optimal or nearly optimal solution to \eqref{dualsdpmarginal}$,$ assume $\U^*_i(x):=\{u\in \U\colon f_i(x,u)=F_i(x)\}, i=1,\ldots,l$ is singleton for $\varphi$-almost all $x\in X,$ let $\rho_i:X\to \U$ be the $\varphi$-measurable mapping such that $\rho_i(x)=\U^*_i(x),$ then
	\begin{equation*}
		z^{\alpha\beta}_{(d^i)^\prime}\to\int_{x\in X}x^\alpha \rho_i(x)^\beta d\varphi(x),~\forall \alpha,\beta\in\N^n_{2(d^i)^\prime},~|\alpha|+|\beta|\le s(2(d^i)^\prime)
	\end{equation*}
	as $(d^i)^\prime\to\infty$.
\end{theorem}

Denoted by $e_k$ the unit vector in $\R^p$ with the $k$-th component equaling to $1.$ According to Theorem \ref{momentconvergence}, it is worth noting that the asymptotic performance of  $u_{(d^i)^\prime}:=\left(z^{0e_1}_{(d^i)^\prime},\ldots,z^{0e_p}_{(d^i)^\prime}\right)$ is meaningless, in other words, the limit of $u_{(d^i)^\prime}$ is not the uncertainty parameter which achieves the optimal.

\begin{remark}{\rm
Fix $d^i$ as the degree of the approximation polynomial for $F_i(x)$. 
For relaxation order $(d^i)^\prime$, if the optimizer $(z_i)^*$ to \eqref{dualsdpmarginal} is flat truncated, then the corresponding $(\mu^i)^*$, namely the solution to the dual problem \eqref{sdpmarginalf} of \eqref{dualsdpmarginal}, constitutes the best approximated $d^i$-degree polynomial $F^{d^i}(x):=\sum\limits_{\alpha\in\N^n_{d^i}}(\mu^i_\alpha)^* x^\alpha\in\R_{d^i}[x]$. 
Besides, under Archimedean condition, the flat truncation for the optimizer to \eqref{dualsdpmarginal} is satisfied generically when the relaxation order is sufficiently large; see, e.g., \cite{nie2014optimality}.
}\end{remark}

Generally, for any $d^i\ge1$ we are able to obtain the best approximated $d^i$-degree polynomial for $F_i(x)$ theoretically. From a computational viewpoint, due to the limits of SDP solvers, we always obtain a suboptimally approximated $d^i$-degree polynomial for $F_i(x)$ with a loose relaxation.

\begin{remark}{\rm 
	Analogical results also hold for the dual problem of \eqref{sdpmarginalg}.
}\end{remark}

\section{The utopia point approximation scalarization for \eqref{RVOP}}\label{sect::4}

Recall the concept of a utopia point given in Definition~\ref{utopia}. 
For $i=1,\ldots,l$, if we approximate $\sup\limits_{u\in \U}f_i(x,u)$ directly with ``joint+marginal" method, we will obtain an upper bound polynomial of $\sup\limits_{u\in \U}f_i(x,u),$ but it is difficult to determine how large $\epsilon$ should be in order to generate a utopia point. 
For this reason, we interchange the order of min-max in \eqref{ideal} and it always holds that
\begin{equation}\label{maxmin}
	\sup_{u\in \U}\min_{x\in \Omega}~f_i(x,u)\le\min_{x\in \Omega}\sup_{u\in \U}~f_i(x,u),
\end{equation}
which implies that the left side of \eqref{maxmin} is literally a utopia point. Moreover, it is reasonable to replace the constraint set $\Omega$ with the bigger and more easily computed set $X$, where $X$ is defined as in \eqref{simpleset}, and it holds the relationship that
\begin{equation}\label{maxminX}
	\sup_{u\in \U}\min_{x\in X}~f_i(x,u)\le\sup_{u\in \U}\min_{x\in \Omega}~f_i(x,u)\le\min_{x\in \Omega}\sup_{u\in \U}~f_i(x,u).
\end{equation}
We then adopt the algorithms proposed in \cite{lasserre2011min} with a loose relaxation to obtain a strictly lower bound of the left side of \eqref{maxminX}, which naturally is a utopia point of the problem~\eqref{RVOP}.

Denoted by $\X_{RP}$ (resp., $\X_{PRP}$) and $\Y_{RP}$ (resp., $\Y_{PRP}$) the sets of all (resp., properly) robust Pareto solutions and (resp., properly) values of the problem~\eqref{UVOP}.
To proceed, we also define the sets of best approximation with utopia point $y^U$ in terms of $\lambda$, $\|\cdot\|_p$ and $F(x):=\left(F_1(x),\ldots,F_l(x)\right)$ given as in the problem~\eqref{VP} by
\begin{align*}
	\mathcal{S}(\lambda,p,F, \Omega):=&\ \left\{x\in\R^n\colon x\in\argmin\{\|\lambda\odot(F(x)-y^U)\|^p_p\colon x\in\Omega\}\right\}, \\
	\mathcal{S}(F,\Omega):=&\ \bigcup_{\lambda\in \ri\Lambda}\bigcup_{1\le p<\infty}\mathcal{S}(\lambda,p,F,\Omega).
\end{align*}
Correspondingly, we define the sets in the objective space for the problem~\eqref{VP},
\begin{align*}
	\A(\lambda,p,F,\Omega):=&\ \left\{F(x)\in\R^l\colon x \in\mathcal{S}(\lambda,p,F,\Omega)\right\}, \\
	\A(F,\Omega):=&\ \bigcup_{\lambda\in \ri\Lambda}\bigcup_{1\le p<\infty}\A(\lambda,p,F,\Omega).
\end{align*}

In what follows, we first give some results on the nonemptiness of $\A(\lambda,p,F,\Omega)$, which coincides with the $\R^l_\geqq$-closedness that is requisite in the proof of Theorem \ref{Utopia} given in \cite[Theorem 4.25]{Ehrgott2005}.
\begin{lemma}\label{lsc}
	For $i=1,\ldots,l$, $F_i(x),$ defined as in \eqref{paraf}$,$ is lower semicontinuous on $X$.
\end{lemma}
\begin{proof}
	For $i=1,\ldots,l$, $f_i(x,u)$ is polynomial (thus continuous) on $X$ for all $u\in \U,$ so the epigraph of which ${\rm epi}(f^u_i):=\left\{(x,t)\in\R^n\times\R\colon t\ge f_i(x,u),~x\in X\right\}$ with respect to $x$ is closed for all $u\in \U$. 
On the other hand, as the epigraph of $F_i(x)$ satisfies ${\rm epi}(F_i)=\bigcap\limits_{u\in \U}{\rm epi}(f^u_i)$, which is the intersection of a family of closed sets, therefore closed, therefore $F_i(x)$ is lower semicontinuous on $X$.
\end{proof}

For fixed $\lambda\in \ri\Lambda$ and $p\in\left[1,\infty\right)$ we define the following function
\begin{align}\label{func4.3}
	\Gamma_{\lambda,p}\colon &X\to \R\\
	&x\mapsto\sum_{i=1}^{l}\left(\lambda_i(F_i(x)-y^U_i)\right)^p. \notag
\end{align}
\begin{lemma}
It holds that $\A(\lambda,p,F,\Omega)\ne\emptyset$ for all $\lambda\in \ri\Lambda$ and $p\in\left[1,\infty\right)$.
\end{lemma}

\begin{proof}
	Notice that the utopia point $y^U < F(x)$ for all $x\in X$, then it holds that
	\begin{equation*}
		\mathcal{S}(\lambda,p,F,\Omega)=\bigl\{x\in\R^n\colon x\in\argmin\{\Gamma_{\lambda,p}(x)\colon x\in\Omega\}\bigr\}.
	\end{equation*}
	For any fixed $\hat{x}\in\Omega$ and any sequence $\{x_k\}_{k\ge1}\subseteq\Omega$ for which $x_k\to\hat{x}$ as $k\to\infty$, we have
	\begin{align*}
		\Gamma_{\lambda,p}(\hat{x})=&\ \sum_{i=1}^{l}\left(\lambda_i(F_i(\hat x)-y^U_i)\right)^p\\
		\le&\ \sum_{i=1}^{l}\left(\lambda_i\left(\liminf_{k\to\infty}F_i(x_k)-y^U_i\right)\right)^p \tag{by Lemma~\ref{lsc}}\\
		\le&\ \liminf_{k\to\infty}\Gamma_{\lambda,p}(x_k),
	\end{align*}
which means $\Gamma_{\lambda,p}$ is lower semicontinuous on $\Omega$. 
On the other hand, $\Omega=\bigcap\limits_{v\in \V}\{x\in\R^n\colon g_j(x,v) \geq0, ~ j=1,\ldots,m\}$, namely the intersection of a family of compact sets in accordance with assumption \textbf{(H1)}, which is compact. Therefore $\mathcal{S}(\lambda,p,F,\Omega)\ne\emptyset$ then $\A(\lambda,p,F,\Omega)\ne\emptyset$.
\end{proof}

\begin{theorem}\label{Utopia}{\rm (c.f. \cite[Theorem 4.25]{Ehrgott2005})}
	\begin{equation*}
		\A(F,\Omega)\subseteq\Y_{PRP}\subseteq\Y_{RP}\subseteq \cl(\A(F,\Omega)).
	\end{equation*}
\end{theorem}

\begin{corollary}\label{Utopiapoints}
	$\mathcal{S}(F,\Omega)\subseteq\X_{PRP}\subseteq\X_{RP}$.
\end{corollary}

Denoted by $F^{d^i}(x), i=1,\ldots,l$, the polynomials obtained from \eqref{sdpmarginalf}, let $a:=\min\{d^1,\ldots,d^l\}$. With a bit abuse of notation, we define the vector polynomial as $\bar{F}^a(x):=(F^{d^1}(x),\ldots,F^{d^l}(x))$. Similarly we define the best approximation sets for $\bar{F}^a(x)$ on $\Omega$:
\begin{align*}
\mathcal{S}(\lambda,p,\bar{F}^a,\Omega):=&\ \left\{x\in\R^n\colon x\in\argmin_{x\in\Omega}\left\{\|\lambda\odot(\bar{F}^a(x)-y^U)\|^p_p \right\}\right\},\\
\mathcal{S}(\bar{F}^a,\Omega):=&\ \bigcup_{\lambda\in \ri\Lambda}\bigcup_{1\le p<\infty}\mathcal{S}(\lambda,p,\bar{F}^a,\Omega),\\
	\A(\lambda,p,\bar F^a,\Omega):=&\ \left\{\bar{F}^a(x)\in\R^l\colon  x\in\mathcal{S}(\lambda,p,\bar{F}^a,\Omega)\right\},\\
	\A(\bar F^a,\Omega):=&\ \bigcup_{\lambda\in \ri \Lambda}\bigcup_{1\le p<\infty}\A(\lambda,p,\bar F^a,\Omega),
\end{align*}
and the function on $X$ defined in terms of $\bar{F}^a,$
\begin{equation*}
	\bar{\Gamma}^a_{\lambda,p}\colon x\mapsto\sum_{i=1}^{l}\left(\lambda_i(\bar F^a_i(x)-y^U_i)\right)^p.
\end{equation*}
Note that $\bar{\Gamma}^a_{\lambda,p}(x)\ge\Gamma_{\lambda,p}(x)$ for all $x\in X$ since $F^{d^i}(x)$ is a upper bound of $F_i(x)$ for $i=1,\ldots,l$.

\begin{lemma}\label{continuos}
	The function $\Gamma_{\lambda,p}(x)$ defined as in \eqref{func4.3} is continuous on $X$ for all $\lambda\in \ri\Lambda$ and $p\in\left[1,\infty\right)$.
\end{lemma}

\begin{proof}
	For each $i=1,\ldots,l,$ it is sufficient to prove $F_i(x)$ is continuous on $X$, and we only need to prove $F_i(x)$ is upper semicontinuous (u.s.c.) on $X$ due to Lemma~\ref{lsc}. 
	Indeed, for all $x\in X$ and all $\{x_k\}\subseteq X$ satisfying $x_k\to x$, there exists $u_k\in \U$ such that
	\begin{equation*}
	F_i(x_k)=f_i(x_k,u_k),
	\end{equation*}
	since $\U$ is compact and $f$ is continuous on $X \times \U.$ 
	Meanwhile, there exists $u\in \U$ which is the accumulation point of $\{u_k\}$ such that
	\begin{equation*}
	\limsup_{k\to\infty}F_i(x_k)=\limsup_{k\to\infty}f_i(x_k,u_k)=f_i(x,u)\le F_i(x),
	\end{equation*}
	therefore $F_i(x)$ is u.s.c. on $X$ and then $\Gamma_{\lambda,p}(x)$ is continuous on $X$.
\end{proof}

Now, by replacing the supremum vector-valued function $F(x)$ with the polynomial approximation $\bar{F}^a(x)$ obtained by the ``joint+marginal" method, we arrive at the following two theorems, which state the convergence conclusions for the utopia point approximation scalarization method with $\bar{F}^a(x)$ on $\Omega$ about the optimal value and optimal points in the sense of (approximated) Pareto, respectively, as $a\to\infty$.
Before that, we also need the following assumption on the robust feasible set $\Omega$ defined as in \eqref{ro_fe_so_se}.
\begin{framed}
	\begin{itemize}
		\item [\textbf{(H2)}] $\Omega$ is the closure of an open set, namely $\cl\bigr({\rm int}(\Omega)\bigl)=\Omega$.
	\end{itemize}
\end{framed}
\begin{remark}{\rm 
	Analogical to the same assumptions presented in \cite{lasserre2012algorithm, jeyakumar2016convergent}, the assumption \textbf{(H2)} plays a significant role in the proofs of Theorems \ref{convergevalue} and \ref{convergenceconstriant}.
}\end{remark}

\begin{theorem}\label{convergevalue}
	Let $\Gamma^*_{\lambda,p}:=\min\limits_{x\in\Omega}\Gamma_{\lambda,p}(x)$ and $(\bar{\Gamma}^t_{\lambda,p})^*:=\min\limits_{x\in\Omega}\bar{\Gamma}_{\lambda,p}(x)$. 
Then for all  $\lambda\in \ri\Lambda$ and $p\in\left[1,\infty\right),$ it holds$:$
	\begin{equation*}
	\liminf\limits_{t\to\infty}(\bar{\Gamma}^t_{\lambda,p})^*=\Gamma^*_{\lambda,p}.
	\end{equation*}
\end{theorem}

\begin{proof}
	With assumption \textbf{(H2)}, for all sufficient small $\alpha > 0,$ the set
	\begin{equation*}
	B^{\lambda,p}_\alpha:=\left\{x\in \Omega\colon \Gamma_{\lambda,p}(x)<\Gamma^*_{\lambda,p}+\alpha\right\}
	\end{equation*}
	is open since $\Gamma_{\lambda,p}(x)$ is continuous on $X$ thanks to Lemma \ref{continuos}, thus is $\varphi$-measurable and $\varphi\bigr(B^{\lambda,p}_\alpha\bigl)=:\kappa>0$.
	
For the sake of readability, we abuse some notation in this proof. 
According to Theorem \ref{convergence}, we know $||F^t_i-F_i||_{L_1(X,\varphi)}\to0$ as $t \to \infty$, for $i=1,\ldots,l$. 
Since $\varphi(\Omega)\le\varphi(X)=1<\infty$, according to Egorov's theorem $F^a_1$ converges to $F_1$ almost uniformly, then there exists a subsequence $\{F^{a(1)}_1\}\subseteq\{F^a_1\}$ such that for the open set $B_1\subseteq \Omega$ satisfying $\varphi(B_1)<\frac{\kappa}{2l}$, $\{F^{a(1)}_1\}$ converges to $F_1$ uniformly on $\Omega\setminus B_1$. 
Similarly, there exists a subsequence $\{F^{a(2)}_2\}\subseteq\{F^{a(1)}_2\}$ such that for the open set $B_2\subseteq \Omega$ satisfying $\varphi(B_2)<\frac{\kappa}{2l}$, $\{F^{a(2)}_2\}$ converges to $F_2$ uniformly on $\Omega\setminus B_2$. 
And so on, denoted by $B:=\bigcup\limits_{i=1}^lB_i$, we obtain $l$ subsequences $\{F^{a(l)}_i\}_{a(l)\ge1}\subseteq\{F^a_i\}_{a\ge1}$ such that $\{F^{a(l)}_i\}$ converges to $F^a_i$ uniformly on the compact set $\Omega\setminus B$, for $i=1,\ldots,l$. 
Therefore there exists $M_i\ge0$ such that $(F^{a(l)}_i(x)-y^U_i)-(F_i(x)-y^U_i)\in[0,M_i]$.
	
Meanwhile for $i=1,\ldots,l$, there exists $0\le a_i\le b_i$ such that $F_i(x)-y^U_i\in[a_i,b_i]$ since $F^{a(l)}_i$ and $F_i$ are continuous on $\Omega\setminus B$. 
Thus
	\begin{align*}
	\bar{\Gamma}^{a(l)}_{\lambda,p}(x)-\Gamma_{\lambda,p}(x)=&\ \sum_{i=1}^{l}\lambda^p_i\left((F^{a(l)}_i(x)-y^U_i)^p-(F_i(x)-y^U_i)^p\right)\\
	\le&\ \sum_{i=1}^{l}\lambda^p_ip(b_i+M_i)^{p-1}\bigl(F^{a(l)}(x)-F_i(x)\bigr),
	\end{align*}
	which means $\{\bar{\Gamma}^{a(l)}_{\lambda,p}(x)\}_{a(l)\ge1}$ converges uniformly to $\Gamma_{\lambda,p}(x)$ on $\Omega\setminus B$.
	
	The open set $B\subseteq \Omega$ satisfies $\varphi(B)\le\sum_{i=1}^{l}\varphi(B_i)<\frac{\kappa}{2}$, therefore
	\begin{align*}
	\varphi(\Omega\setminus B)>\varphi(\Omega)-\frac{\kappa}{2}.
	\end{align*}
	Then $B^{\lambda,p}_\alpha \cap \big(\Omega\setminus B\big)\ne\emptyset$, and there exists $x\in B^{\lambda,p}_\alpha$ such that $\lim\limits_{t(l)\to\infty}\bar{\Gamma}^{t(l)}_{\lambda,p}(x)=\Gamma_{\lambda,p}(x)<\Gamma^*_{\lambda,p}+\alpha$, which means
	\begin{equation*}
	\liminf\limits_{t\to\infty}(\bar{\Gamma}^t_{\lambda,p})^*\le\lim\limits_{t(l)\to\infty}\bar{\Gamma}^{t(l)}_{\lambda,p}(x)<\Gamma^*_{\lambda,p}+\alpha.
	\end{equation*}
	As $\alpha>0$ is arbitrary, we have $\liminf\limits_{t\to\infty}(\bar{\Gamma}^t_{\lambda,p})^*\le\Gamma^*_{\lambda,p}$. On the other hand, $\liminf\limits_{t\to\infty}(\bar{\Gamma}^t_{\lambda,p})^*\ge\Gamma^*_{\lambda,p}$ for $\bar{\Gamma}_{\lambda,p}(x)\ge\Gamma_{\lambda,p}(x)$ for all $x\in X$, we thus complete the proof.
\end{proof}

\begin{remark}{\rm 
We would mention here again that, the assumption \textbf{(H2)} is necessary for Theorem \ref{convergevalue}, 
and without it, the conclusion may not hold. 
For example,
let $X=[0,3]\subset\R$ and $\Omega=:[1,2]\cap\{3\}\subset X$, then $[1,2]=\cl\big({\rm int}(\Omega)\big)\ne\Omega$ thus assumption \textbf{(H2)} is no longer satisfied. 
Without loss of generality, let $F(x)=-x$, then $F(x)$ admits its minimum on $\Omega$ at $x^*=3$, while for $0<\alpha<1$, $B_\alpha:=\{x\in\Omega\colon F(x)<F(x^*)+\alpha\}=\{3\}$ is not open.
}\end{remark}

We have to point out that, it is usually impractical for us to compute $\liminf\limits_{a\to\infty}(\bar{\Gamma}^a_{\lambda,p})^*$ directly due to the fact that we never have information about the future. 
Instead we compute $\lim\limits_{k\to\infty}\inf\limits_{1\le a\le k}\{(\bar{\Gamma}^a_{\lambda,p})^*\}$, and because of
\begin{equation*}
	\Gamma^*_{\lambda,p}\le\lim\limits_{k\to\infty}\inf\limits_{1\le a\le k}\{(\bar{\Gamma}^a_{\lambda,p})^*\}\le\liminf\limits_{a\to\infty}(\bar{\Gamma}^a_{\lambda,p})^*=\Gamma^*_{\lambda,p},
\end{equation*}
we actually obtain the same minimum value. With a bit abuse of notation, define
\begin{equation*}
	(\bar{\Gamma}^{a(k)}_{\lambda,p})^*:=\inf\limits_{1\le a\le k}\{(\bar{\Gamma}^a_{\lambda,p})^*\},
\end{equation*}
we therefore have the convergence conclusion for optimal points in the sense of Pareto as below.

\begin{theorem}\label{pointconvergence1}
	Obeying the notation defined above for $a(k),$ it holds that
	\begin{equation*}
		\emptyset\ne\bigcap\limits_{r=1}^{\infty}\cl\bigl(\bigcup\limits_{a(k)=r}^{\infty}\mathcal{S}(\lambda,p,\bar F^{a(k)},\Omega)\bigr)\subseteq\mathcal{S}(\lambda,p,F,\Omega).
	\end{equation*}
\end{theorem}

\begin{proof}
We firstly prove $\bigcap\limits_{r=1}^{\infty}\cl\bigl(\bigcup\limits_{a(k)=r}^{\infty}\mathcal{S}(\lambda,p,\bar F^{a(k)},\Omega)\bigr)\ne\emptyset$. 
Indeed, since $\bar{\Gamma}^{t(k)}_{\lambda,p}(x)$ is continuous on the compact set $\Omega$, it is obvious that $\mathcal{S}(\lambda,p,\bar F^{a(k)},\Omega)\ne\emptyset$, then $\cl\bigl(\bigcup\limits_{{t(k)}=r}^{\infty}\mathcal{S}(\lambda,p,\bar F^{a(k)},\Omega)\bigr)\ne\emptyset$ and $\cl\bigl(\bigcup\limits_{{a(k)}=r}^{\infty}\mathcal{S}(\lambda,p,\bar F^{a(k)},\Omega)\bigr)\supseteq \cl\bigl(\bigcup\limits_{{a(k)}=r+1}^{\infty}\mathcal{S}(\lambda,p,\bar F^{a(k)},\Omega)\bigr)$, therefore there is at least one point in the intersection of a sequence of monotonically decreasing compact sets.

Now, we are going to show $\bigcap\limits_{r=1}^{\infty}\cl\bigl(\bigcup\limits_{a(k)=r}^{\infty}\mathcal{S}(\lambda,p,\bar F^{a(k)},\Omega)\bigr)\subseteq\mathcal{S}(\lambda,p,F,\Omega)$. 
Indeed, for all $\hat{x}\in\bigcap\limits_{r=1}^{\infty}\cl\bigl(\bigcup\limits_{a(k)=r}^{\infty}\mathcal{S}(\lambda,p,\bar F^{a(k)},\Omega)\bigr)$, there exists a sequence $\{x_{{a(k)}^\prime}\}$ satisfying $x_{{a(k)}^\prime}\in\mathcal{S}(\lambda,p,\bar F^{{a(k)}^\prime},\Omega)$ and they accumulate to $\hat{x}$. Then for ${a(k)}^\prime\ge1$, it holds
	\begin{equation*}
	\Gamma^*_{\lambda,p}\le\Gamma_{\lambda,p}(x_{{a(k)}^\prime})\le\bar\Gamma^{{a(k)}^\prime}_{\lambda,p}(x_{{t(k)}^\prime})=(\bar{\Gamma}^{{a(k)}^{\prime}}_{\lambda,p})^*.
	\end{equation*}
	Applying the continuity of $\Gamma_{\lambda,p}(x)$ and Theorem~\ref{convergevalue}, we prove the conclusion.
\end{proof}

\begin{remark}{\rm
	Indeed, the upper limit closure set in Theorem \ref{pointconvergence1} is literally the set of all accumulation points for $\{\mathcal{S}(\lambda,p,\bar F^{a(k)},\Omega)\}_k$ as $k\to \infty.$
	Theorem \ref{pointconvergence1} tells us, in other words, that as we increasing the order of approximation polynomials, the limiting point of $\mathcal{S}(\lambda,p,\bar F^{a(k)},\Omega)$ always exists and belongs to $\mathcal{S}(\lambda,p,F,\Omega).$
}\end{remark}

Denote by $G^{e^j}(x), j=1,\ldots,m$ the polynomials obtained from \eqref{sdpmarginalg}, let $b:=\{e^1,\ldots,e^m\}$.
Again, with a bit abuse of notation, we define the set 
\[
\bar{\Omega}^b:=\left\{x\in\R^n\colon G^{e^j}(x)\ge0,~j=1,\ldots,m\right\}\cap X.
\]
Then we have the following vector polynomial optimization problem,
\begin{align*}\label{VOPab}
\min~&\bar{F}^a(x)\tag{VOP$^{a,b}$}\\
{\rm s.t.}~&x\in\bar{\Omega}^b.
\end{align*}
Analogically we define the best approximation sets for the problem \eqref{VOPab}
\begin{align*}
	\mathcal{S}(\lambda,p,\bar{F}^a,\bar{\Omega}^b):=&\ \bigl\{x\in\R^n\colon x\in\argmin\{||\lambda\odot(\bar{F}^a(x)-y^U)||_p\colon x\in\bar{\Omega}^b\}\bigr\},\\
	\mathcal{S}(\bar{F}^a,\bar{\Omega}^b):=&\ \bigcup_{\lambda\in \ri \Lambda}\bigcup_{1\le p<\infty}\mathcal{S}(\lambda,p,\bar{F}^a,\bar{\Omega}^b),\\
	\A(\lambda,p,\bar F^a,\bar{\Omega}^b):=&\ \left\{\bar{F}^a(x)\in\R^l\colon  x\in\mathcal{S}(\lambda,p,\bar{F}^a,\bar{\Omega}^b)\right\},\\
	\A(\bar F^a,\bar{\Omega}^b):=&\ \bigcup_{\lambda\in \ri \Lambda}\bigcup_{1\le p<\infty}\A(\lambda,p,\bar F^a,\bar{\Omega}^b).
\end{align*}
Note that $\bar{\Omega}^b\subseteq\Omega$ since $G^{e^j}(x)$ is a lower bound of $G^j(x)$ (defined as in \eqref{parag}) for $j=1,\ldots,m$.

In what follows, we give the convergence conclusions for the problem~\eqref{VOPab} as $a,b\to\infty$.

\begin{theorem}\label{convergenceconstriant}
Assume that ${\rm int}\big(\{x\in\R^n\colon G_j(x)=0\}\big)\cap X=\emptyset$ for all $j=1,\ldots,m,$ then it holds
	\begin{equation*}
		\cl\bigl(\lim_{b\to\infty}\bar{\Omega}^b\bigr)=\Omega.
	\end{equation*} 
\end{theorem}
 
\begin{proof}
By assumption \textbf{(H1)} we know that $\Omega$ is compact, all we need (thanks to assumption \textbf{(H2)}) is to prove 
\[
{\rm int}(\Omega)\subseteq \cl\left(\bigcup\limits_{r=1}^{\infty}\bigcap\limits_{b=r}^{\infty}\bar{\Omega}^b\right)\subseteq \cl\bigl(\lim\limits_{b\to\infty}\bar{\Omega}^b\bigr) \subseteq\bigcap\limits_{r=1}^{\infty}\cl\bigl(\bigcup\limits_{b=r}^{\infty}\bar{\Omega}^b\bigr)\subseteq\Omega.
\] 
By contradiction, suppose there exists $x_0\in {\rm int}(\Omega)$ such that 
$x_0 \notin \cl\left(\bigcup\limits_{r=1}^{\infty}\bigcap\limits_{b=r}^{\infty}\bar{\Omega}^b\right)$, 
then $\exists \alpha > 0$ such that 
$B_\alpha(x_0):=\left\{x\in\R^n\colon\| x-x_0\| < \alpha\right\}\subseteq {\rm int}(\Omega)$ 
but $B_\alpha(x_0)\cap\left(\cl\bigl(\bigcup\limits_{r=1}^{\infty}\bigcap\limits_{b=r}^{\infty}\bar{\Omega}^b\bigr)\right)=\emptyset$, which means that there exists a sequence $\{b_k\}\subseteq\N^+$ such that $B_\alpha(x_0)\cap\bar{\Omega}^{b_k}=\emptyset$. 
Furthermore, for some $j\in\{1,\ldots,m\}$, there exists an open set $O_\alpha(x_0)\subseteq B_\alpha(x_0)$ such that $G_j(x)>0$ and $G^{e^j}(x)<0$ for all $x\in O_\alpha(x_0)$, where $\bar{\Omega}^{b_k}=\left\{x\in\R^n\colon G^{e^j}(x)\ge0,~j=1,\ldots,m\right\}\cap X$. 
Since $\varphi(O_\alpha(x_0)) > 0$, this contradicts to the fact that $\|G^{e^j}-G_j\|_{L_1(X,\varphi)}\to0$ given in Theorem \ref{convergence}, thus ${\rm int}(\Omega)\subseteq \cl\left(\bigcup\limits_{r=1}^{\infty}\bigcap\limits_{b=r}^{\infty}\bar{\Omega}^b\right)$. 

On the other hand, it is obvious that 
$\bigcap\limits_{r=1}^{\infty}\cl\bigl(\bigcup\limits_{b=r}^{\infty}\bar{\Omega}^b\bigr)\subseteq\Omega$ 
since $\bar{\Omega}^b\subseteq\Omega$. 
Therefore $\cl\bigl(\lim\limits_{b\to\infty}\bar{\Omega}^b\bigr)=\Omega$ as desired.
\end{proof}

\begin{theorem}\label{VOPabapprox}
	With the same assumption as in {\rm Theorem~\ref{convergenceconstriant},} we can solve \eqref{RVOP} approximately via \eqref{VOPab}$,$ mathematically it holds$:$
	\begin{equation*}
	\bigcap\limits_{r^\prime=1}^\infty \cl\biggr(\bigcup\limits_{b=r^\prime}^\infty\bigcap\limits_{r=1}^{\infty}\cl\bigl(\bigcup\limits_{a(k)=r}^{\infty}\mathcal{S}(\lambda,p,\bar F^{a(k)},\bar\Omega^b)\bigr)\biggl)=
\bigcap\limits_{r=1}^{\infty}\cl\bigg(\bigcup\limits_{a(k)=r}^{\infty}\mathcal{S}(\lambda,p,\bar F^{a(k)},\Omega)\bigg).
	\end{equation*}
\end{theorem}

\begin{proof}
First of all, $\bar{\Omega}^b$ is compact since $G^{b(k)}$ is continuous and $\bar{\Omega}^b\subseteq\Omega$, thus $\mathcal{S}(\lambda,p,\bar F^{a(k)},\bar\Omega^b)\ne\emptyset$ and $\bigcap\limits_{r=1}^{\infty}\cl\bigl(\bigcup\limits_{a(k)=r}^{\infty}\mathcal{S}(\lambda,p,\bar F^{a(k)},\bar\Omega^b)\bigr)\ne\emptyset$ for all $b$.
	
Let $\hat{x}\in\bigcap\limits_{r=1}^{\infty}\cl\bigg(\bigcup\limits_{a(k)=r}^{\infty}\mathcal{S}(\lambda,p,\bar F^{a(k)},\Omega)\bigg)$ be any given, 
for all $\alpha>0$, there always exists a sequence $\left\{\mathcal{S}(\lambda,p,\bar F^{a(k)},\Omega) \right\}$ such that 
$\left\{x\in\R^n\colon \|x-\hat{x}\|<\alpha \right\}=:B_\alpha(\hat{x})\cap\mathcal{S}(\lambda,p,\bar F^{a(k)},\Omega)\ne\emptyset.$ 
Meanwhile, for $\alpha$ defined above there also exists a sequence $\{\bar{\Omega}^b\}$ such that $B_\alpha(\hat{x})\cap\bar{\Omega}^b\ne\emptyset$ 
according to Theorem~\ref{convergenceconstriant}. 
Since $\bar \Gamma_{\lambda,p}(x)$ is continuous on $\Omega$, when $\alpha$ is sufficiently small there always exists 
$\bar{\Omega}^b$ such that 
$\argmin\limits_{x\in\bar{\Omega}^b}\bar \Gamma_{\lambda,p}(x)\subseteq B_\alpha(\hat x),$
then there is a subsequence $\{a(k)^\prime\}\subseteq\{a(k)\}$ such that 
$\argmin\limits_{x\in\bar{\Omega}^b}\bar{\Gamma}^{a(k)^\prime}_{\lambda,p}(x)\subseteq B_\alpha(\hat x),$ 
namely 
$\hat{x}\in\bigcap\limits_{r^\prime=1}^\infty \cl\biggr(\bigcup\limits_{b=r^\prime}^\infty\bigcap\limits_{r=1}^{\infty}\cl\bigl(\bigcup\limits_{a(k)=r}^{\infty}\mathcal{S}(\lambda,p,\bar F^{a(k)},\bar\Omega^b)\bigr)\biggl),$ 
thus 
\[
\bigcap\limits_{r=1}^{\infty}\cl\left(\bigcup\limits_{a(k)=r}^{\infty}\mathcal{S}(\lambda,p,\bar F^{a(k)},\Omega)\right)\subseteq\bigcap\limits_{r^\prime=1}^\infty \cl\left(\bigcup\limits_{b=r^\prime}^\infty\bigcap\limits_{r=1}^{\infty}\cl\left(\bigcup\limits_{a(k)=r}^{\infty}\mathcal{S}(\lambda,p,\bar F^{a(k)},\bar\Omega^b)\right)\right).
\]
	
Conversely, let $\hat{x}\in\bigcap\limits_{r^\prime=1}^\infty \cl\biggr(\bigcup\limits_{b=r^\prime}^\infty\bigcap\limits_{r=1}^{\infty}\cl\bigl(\bigcup\limits_{a(k)=r}^{\infty}\mathcal{S}(\lambda,p,\bar F^{a(k)},\bar\Omega^b)\bigr)\biggl)$ be any given, 
for all $\alpha>0$, there always exists a sequence $\{b_k\}$ such that 
$B_\alpha(\hat x)\cap\bigg(\bigcap\limits_{r=1}^{\infty}\cl\bigl(\bigcup\limits_{a(k)=r}^{\infty}\mathcal{S}(\lambda,p,\bar F^{a(k)},\bar\Omega^{b_k})\bigr)\bigg)\ne\emptyset$, 
thus there is a subsequence $\{a(k)^\prime\}\subseteq\{a(k)\}$ such that 
$B_\alpha(\hat x)\cap\mathcal{S}(\lambda,p,\bar F^{a(k)^\prime},\bar\Omega^{b_k})\ne\emptyset.$
According to Theorem~\ref{convergenceconstriant},
$B_\alpha(\hat x)\cap\mathcal{S}(\lambda,p,\bar F^{a(k)^\prime},\Omega)\ne\emptyset$ 
since $\bar \Gamma^{a(k)^\prime}_{\lambda,p}$ are continuous on $\Omega$, moreover 
$\hat{x}\in \bigcap\limits_{r=1}^{\infty}\cl\left(\bigcup\limits_{a(k)=r}^{\infty}\mathcal{S}(\lambda,p,\bar F^{a(k)},\Omega)\right),$ 
then 
\[
\bigcap \limits_{r^\prime=1}^\infty \cl \biggr(\bigcup\limits_{b=r^\prime}^\infty\bigcap\limits_{r=1}^{\infty}\cl\bigl(\bigcup\limits_{a(k)=r}^{\infty}\mathcal{S}(\lambda,p,\bar F^{a(k)},\bar\Omega^b)\bigr)\biggl)\subseteq\bigcap\limits_{r=1}^{\infty}\cl\left(\bigcup\limits_{a(k)=r}^{\infty}\mathcal{S}(\lambda,p,\bar F^{a(k)},\Omega)\right).
\] 
The proof is completed.
\end{proof}

We end this section by mentioning that (by Theorem \ref{pointconvergence1}), for all $\lambda\in\Lambda$ and $p\in\left[1,\infty\right),$
\begin{equation*}
\emptyset\ne\bigcap\limits_{r^\prime=1}^\infty \cl\biggr(\bigcup\limits_{b=r^\prime}^\infty\bigcap\limits_{r=1}^{\infty}\cl\bigl(\bigcup\limits_{a(k)=r}^{\infty}\mathcal{S}(\lambda,p,\bar F^{a(k)},\bar\Omega^b)\bigr)\biggl)\subseteq\mathcal{S}(\lambda,p,F,\Omega).
\end{equation*}
Furthermore, by invoking Theorem~\ref{Utopiapoints}, we also have the following relationship,
\begin{equation*}
	\bigcup_{\lambda\in \ri\Lambda}\bigcup_{1\le p<\infty}\bigcap\limits_{r^\prime=1}^\infty \cl\biggr(\bigcup\limits_{b=r^\prime}^\infty\bigcap\limits_{r=1}^{\infty}\cl\bigl(\bigcup\limits_{a(k)=r}^{\infty}\mathcal{S}(\lambda,p,\bar F^{a(k)},\bar\Omega^b)\bigr)\biggl)\subseteq\X_{PRP}\subseteq\X_{RP}.
\end{equation*}
Correspondingly, in the objective space it holds,
\begin{equation*}
	\bigcup_{\lambda\in \ri\Lambda}\bigcup_{1\le p<\infty}\bigcap\limits_{r^\prime=1}^\infty \cl\biggr(\bigcup\limits_{b=r^\prime}^\infty\bigcap\limits_{r=1}^{\infty}\cl\bigl(\bigcup\limits_{a(k)=r}^{\infty}\mathcal{A}(\lambda,p,\bar F^{a(k)},\bar\Omega^b)\bigr)\biggl)\subseteq\Y_{PRP}\subseteq\Y_{RP}.
\end{equation*}

\section{The Scheme and Examples}\label{sect::5}

As we mentioned before, all of the accumulation points of $\{\mathcal{S}(\lambda,p,F^{a(k)},\bar\Omega^b)\}_{k,b\ge1}$ are properly robust Pareto point to the problem~\eqref{UVOP} for all fixed $\lambda$ and $p$. 
From the viewpoint of computation, we could fix some appropriate orders $a$ and $b$ for the approximating polynomials, then solving the problem~\eqref{VOPab}. 
The procedure of our algorithm is given below.
\medskip

\begin{algorithm}[H]
	\caption{The Utopia Point Method-Based Approximation}
	\SetKwInOut{Input}{Input}
	\SetKwInOut{Output}{Output}
	\Input
	{
		Appropriate orders of approximating polynomials: $a$ and $b$\\
		Maximal order of $\ell_p$-norm: $lp$\\
		A discrete set: $\Lambda\subseteq(0,1)$
	}
	\Output{
		Properly robust Pareto value: $pareto$
	}
	\BlankLine
	Set $pareto=\emptyset$\\
	Compute \textit{utopia point} $y^U$ via \eqref{maxminX} with algorithm proposed in \cite{lasserre2011min}\\
	Generate approximation vector polynomials $\bar F^a$ and $\bar G^b$ via SDP \eqref{sdpmarginalf} and \eqref{sdpmarginalg} \cite{cvx,gb08} respectively\\
	\For{$p=1:lp$}{
		\For{$\lambda\in\Lambda$}{
			compute $\mathcal{S}(\lambda,p,\bar F^a,\bar \Omega^b)$ for \eqref{VOPab} with Lasserre's hierarchy \cite{henrion2009gloptipoly,doi:10.1080/10556789908805766}\\
			\For{$x\in\mathcal{S}(\lambda,p,\bar F^a,\bar \Omega^b)$}{
				$pareto\longleftarrow pareto\bigcup\{F(x)\}$
			}
		}
	}
\end{algorithm}
\bigskip
In what follows, we design some examples to illustrate our algorithm.

\begin{example}\rm
	Consider the robust vector polynomial optimization problem
	\begin{align}\label{ex1}
	\min_{x}~&\sup_{u\in \U}f(x,u)\\
	{\rm s.t.}~&\inf_{v\in \V}g(x,v)\ge0\notag.
	\end{align}
	Let $f(x,u):=(f_1(x,u),f_2(x,u))$, $f_i\colon\R^2\times\R\to\R$, $i=1,2$, where
	\begin{equation*}
	f_1(x,u) = x^2_1u^2+x^2_2u, \quad
	f_2(x,u)=-x^4_1-2x^2_1x^2_2-x^4_2+x^2_1-u^2.
	\end{equation*}
	The uncertainty parameter $u\in \U:=\{u\in\R\colon0\le u\le1\}$. 
Let $g\colon\R^2\times\R^2\to\R$ be the constraint function given as
	\begin{equation*}
	g(x,v):=-x^2_1v^2_2-x^2_2v^2_1+1,
	\end{equation*}
	where the uncertainty set $\V:=\{(v_1,v_2)\in\R^2\colon1\le v_1\le\sqrt{2},~1\le v_2\le\sqrt{2}\}$. The compact simple set $X:=\{(x_1,x_2)\in\R^2\colon x^2_1+x^2_2\le2\}$ covers all feasible $x$, and we can easily compute $F_1$ and $F_2$, the supremums of $f_1$ and $f_2$ for all $x\in X$ over $\U$:
	\begin{align*}
	F_1(x):=&\ \sup_{u\in \U}f_1(x,u)=x^2_1+x^2_2,\\
	F_2(x):=&\ \sup_{u\in \U}f_2(x,u)=-x^4_1-2x^2_1x^2_2-x^4_2+x^2_1.
	\end{align*}
	Moreover, the infimum of $g$ for all $x\in X$ over $\V$ is
	\begin{equation*}
		G(x):=\inf_{v\in V}g(x,v)=1-x^2_1-x^2_2.
	\end{equation*}
	We can depict the robust Pareto front of \eqref{ex1} in advance.
	\begin{figure}[htbp]
		\includegraphics*[width=3in]{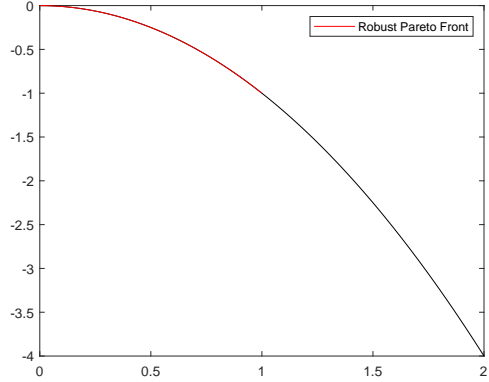}
		\caption{Robust Pareto Front}\label{fig:1}
	\end{figure}
	
	The red curve in {\sc Figure}~\ref{fig:1} is the robust Pareto front of \eqref{ex1} and the black curve represents the southwest lower bound of $(F_1(x),F_2(x))\in\R^2$ for all $x\in X$.
	
	We approximate $F_1(x)$ and $G(x)$ with 2-degree polynomials, and approximate $F_2(x)$ with 4-degree polynomial. With no more than 7-order relaxation, we obtian
	\begin{align*}
		\bar{F}^2_1(x)=&\ \tfrac{2251799814263309}{2251799813685248} x_1^2 + \tfrac{4503599606726343}{4503599627370496} x_2^2 + \tfrac{1501956633657483}{151115727451828646838272}\\
		\approx&\ x^2_1+x^2_2, \\
		\bar{F}^4_2(x)=&\ -\tfrac{70368744161045}{70368744177664} x_1^4 - \tfrac{9007199252613767}{4503599627370496} x_1^2x_2^2 - \tfrac{70368744161045}{70368744177664} x_2^4 + \tfrac{1125899906575295}{1125899906842624} x_1^2  \\
		&\ + \tfrac{9007199250463721}{18014398509481984} x_2^2 + \tfrac{1105826627040601}{4835703278458516698824704}\\
		\approx&\ -x^4_1-2x^2_1x^2_2-x^4_2+x^2_1+0.5x_2^2, \\
		\bar{G}^2(x)=&\ -\tfrac{8896505411494041}{9007199254740992} x_1^2 - \tfrac{4447374274883833}{4503599627370496} x_2^2 - \tfrac{5753780886491375}{21778071482940061661655974875633165533184} x_1  \\ 
		&\ -\tfrac{5752012769693551}{21778071482940061661655974875633165533184} x_2
		 + \tfrac{5980454757129867}{174224571863520493293247799005065324265472} x_1x_2 \\
		 &\ +\tfrac{8794748612256617}{9007199254740992}\\
		\approx&\ -0.9877x^2_1-0.9875x^2_2+0.9764.
	\end{align*}
	
	Eventually, with the approximation polynomials we obtained, we have the following vector polynomial optimization problem
	
	\begin{align}\label{ex1approx}
		\min_{x}~& \left(\bar{F}^2_1(x),\bar{F}^4_2(x) \right)\\
		{\rm s.t.}~&\bar{G}^2(x)\ge0\notag\\
			&x\in X\notag.
	\end{align}
	For $\lambda\in[0.01,0.99]$ and $1\le lp\le12$, we compute \eqref{ex1approx} with utopia point method-based approximation via Lasserre's hierarchy, then we obtian the following result.
	\begin{figure}[htbp]
		\includegraphics*[width=3in]{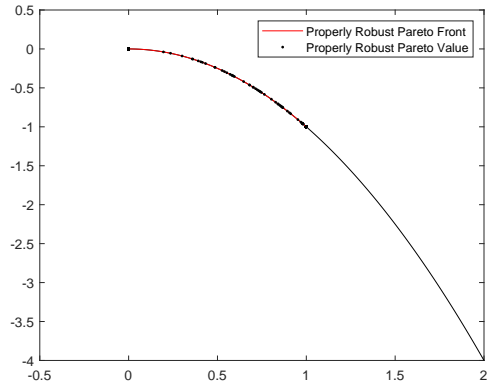}
		\caption{Robust Pareto Front and Properly Robust Pareto Values}\label{fig:2}
	\end{figure}
	
	The red curve in {\sc Figure}~\ref{fig:2} represents the robust Pareto front of \eqref{ex1}, the black dots are properly robust Pareto values we obtained.
\end{example}

Finally, we consider a particular case, i.e., dropping the uncertainty parameters, and the problem~\eqref{UVOP} will be reduced to a vector nonconvex polynomial optimization problem. 
In this case, our algorithm is still available if only omitting the approximation step.

\begin{example}\rm
	Consider the vector polynomial optimization problem
	\begin{align}\label{ex2}
	\min_{x}~&f(x)\\
	{\rm s.t.}~&g(x)\ge0\notag\\
		~&x\in X.\notag
	\end{align}
	Let $f(x):=\left(f_1(x),f_2(x)\right)$, $f_i\colon\R^2\to\R$, $i=1,2$ with
	\begin{align*}
	f_1(x) = &x_1+x_2, \\
	f_2(x) = &-x^2_1+2x_1x_2+0.5x_1+0.5x_2-0.0625.
	\end{align*}
	The constriant set $X:=\{(x_1,x_2)\colon0\le x_1,x_2\le1\}$, and
	\begin{equation*}
		g_1(x)=x^2_1+x^2_2+2x_1x_2-x_1-x_2.
	\end{equation*}
	With the weights $\lambda$ chosen in $[0.01, 0.99]$, we approximate the Pareto front of the problem~\eqref{ex2}; see the following {\sc Figure}~\ref{fig:3456}.
	\begin{figure}[htbp]
		\subfigure[lp=12]{
		\includegraphics*[width=2.5in]{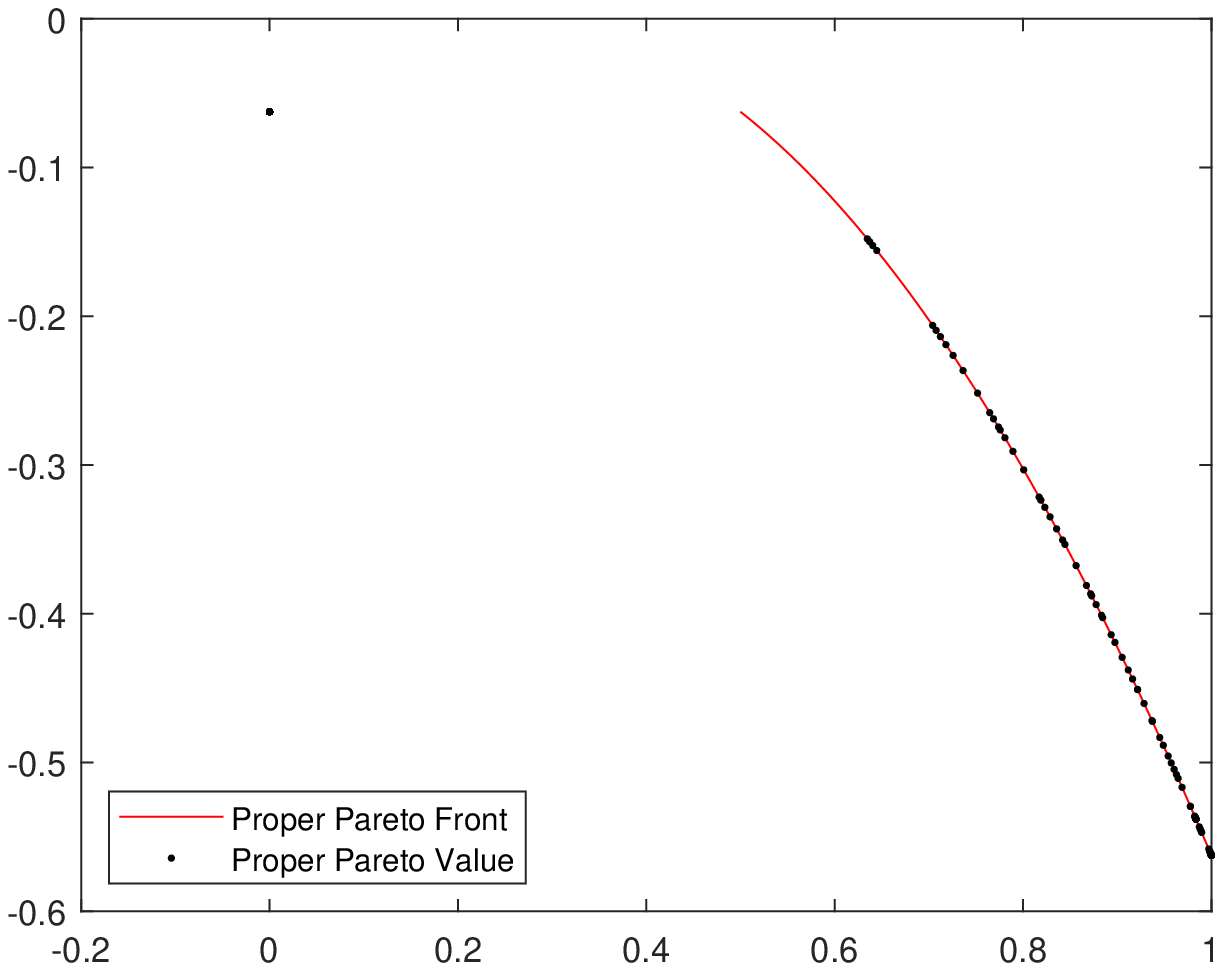}\label{fig:3}
		}
		\quad
		\subfigure[lp=18]{
		\includegraphics*[width=2.5in]{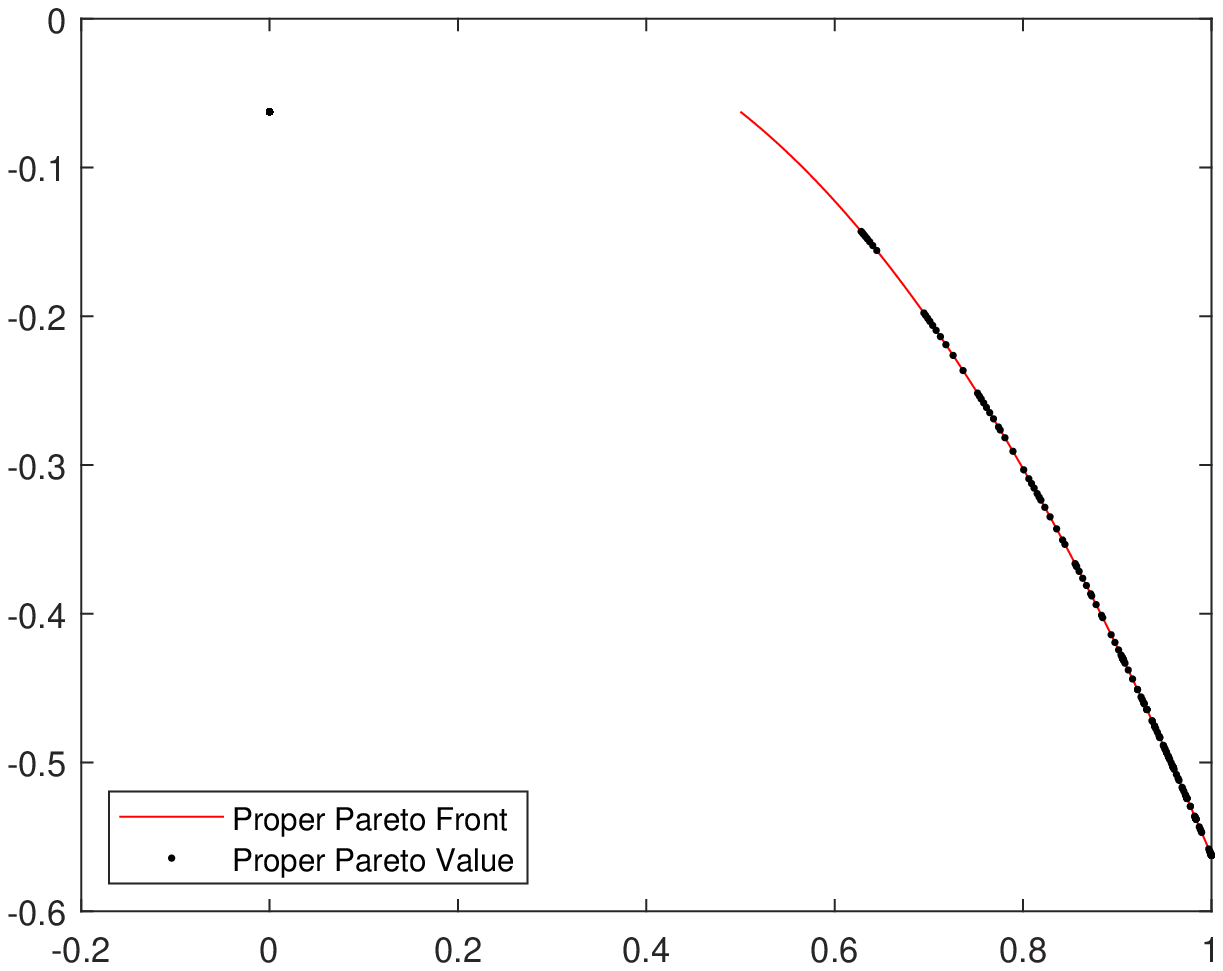}\label{fig:4}
		}
		\quad
		\subfigure[lp=24]{
		\includegraphics*[width=2.5in]{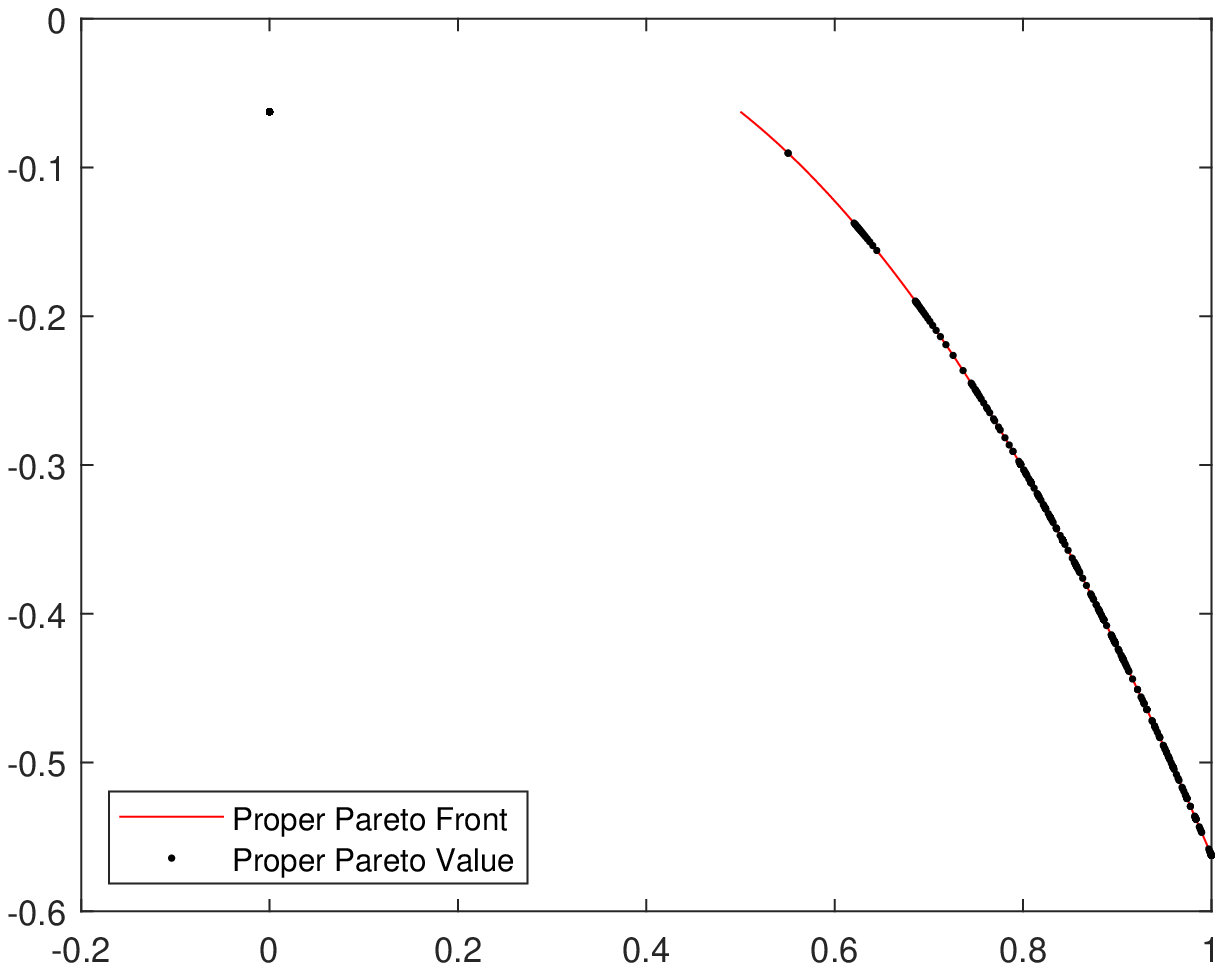}\label{fig:5}
		}
		\quad
		\subfigure[lp=30]{
		\includegraphics*[width=2.5in]{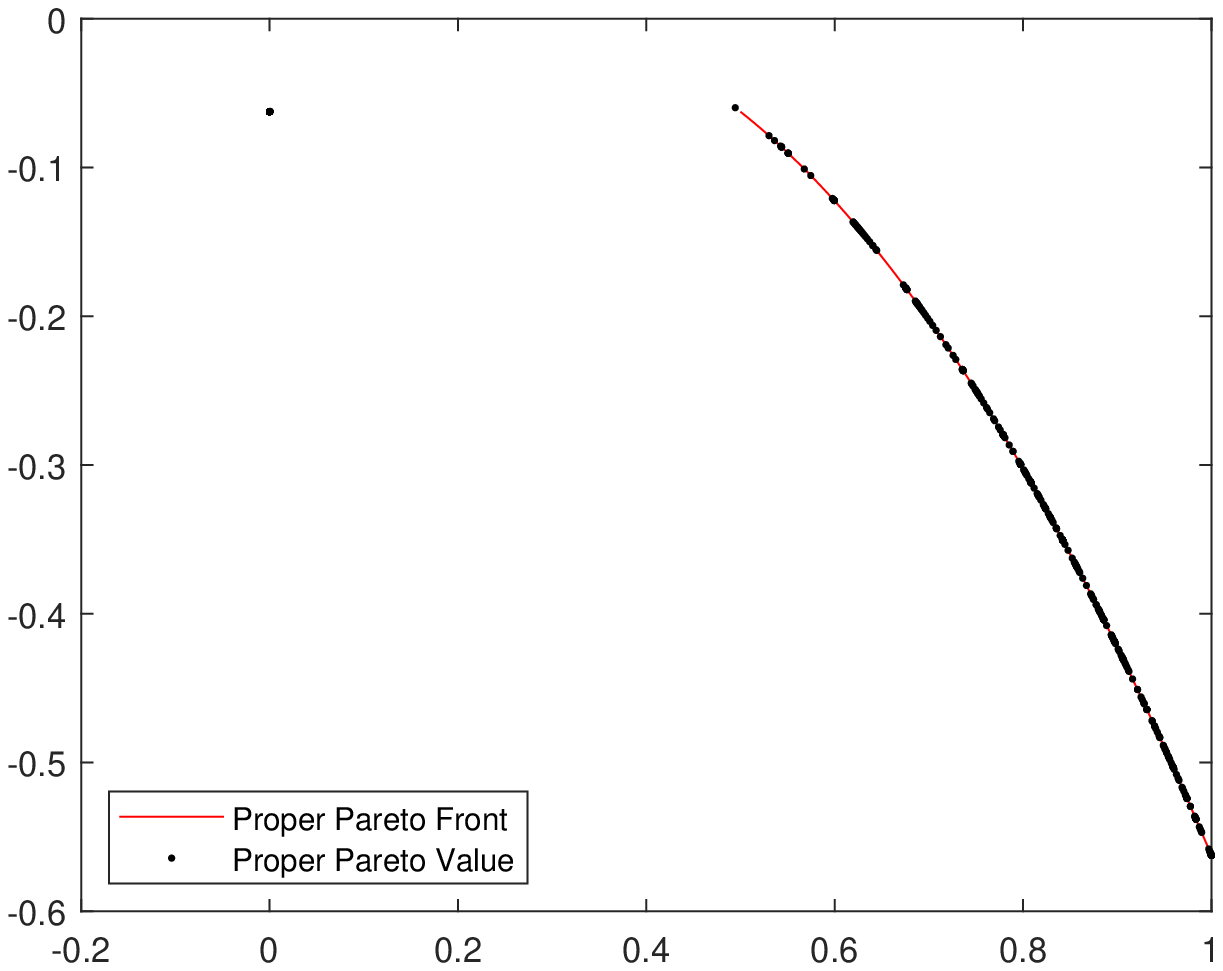}\label{fig:6}
		}
		\caption{Pareto Front and Properly Pareto Values}
		\label{fig:3456}
	\end{figure}
	
	The red curves in {\sc Figure}~\ref{fig:3456} represent the proper Pareto front of the problem~\eqref{ex2} and the black dots are solutions of our algorithm. 
In fact, we are capable to approximate the closure of the proper Pareto front as we increase the order of $p$-norm.
\end{example}

\section{Conclusion}\label{sect::6}
In this paper, we mainly studied a class of robust vector polynomial optimization problems without any convex assumptions. 
We solved such problems by combining/improving the utopia point method and ``joint+marginal" relaxation method.
Particularly, comparing with our previous works \cite{Jiao2019,Jiao2020,Lee2018,Lee2021,Lee2021JOGO}, we developed the previous results from convex case to nonconvex case.

Besides, it would be of great interest to examine in what situation the assumption in Theorem \ref{convergenceconstriant} doesn't hold.
Under the consideration that $\V$ is set to be compact, it seems that the polynomial who breaks the assumption in Theorem~\ref{convergenceconstriant} should have the form $g(x,v)=\sum_{k=1}^{d}\tilde{g}_k(x)\bar{g}_k(v).$
A deeper charaterization will be carried out in the future studies.

\subsection*{Acknowledgments}
This work was supported by the Key Project of National Natural Science Foundation of China under Grant NSFC12031016, Grant 12026607, and the Pazhou lab, GuangZhou; by the National Key R\&D Program of Chia under Grant 2020YFA0712203, Grant 2020YFA0712201.
Liguo Jiao was supported by Natural Science Foundation of Jilin Province (YDZJ202201ZYTS302). 
Jae Hyoung Lee was supported by the National Research Foundation of Korea (NRF)
grant funded by the Korea government (MSIT) (NRF-2021R1C1C2004488).


\begin{thebibliography}{10}

\bibitem{Beck2009}
A.~Beck and A.~Ben-Tal.
\newblock Duality in robust optimization: primal worst equals dual best.
\newblock {\em Operations Research Letters}, 37:1--6, 2009.

\bibitem{Ben-Tal2009}
A.~Ben-Tal, L.~E. Ghaoui, and A.~Nemirovskii.
\newblock {\em Robust Optimization}.
\newblock Princeton University Press, Princeton, NJ, 2009.

\bibitem{Ben-Tal2002}
A.~Ben-Tal and A.~Nemirovski.
\newblock Robust optimization---methodology and applications.
\newblock {\em Mathematical Programming}, 92:453--480, 2002.

\bibitem{Ben-Tal2008}
A.~Ben-Tal and A.~Nemirovski.
\newblock Selected topics in robust convex optimization.
\newblock {\em Mathematical Programming}, 112(1):125--158, 2008.

\bibitem{Bertsimas2011}
D.~Bertsimas, D.~B. Brown, and C.~Caramanis.
\newblock Theory and applications of robust optimization.
\newblock {\em SIAM Review}, 53:464--501, 2011.

\bibitem{Bitran1980}
G.~R. Bitran.
\newblock Linear multiple objective problems with interval coefficients.
\newblock {\em Management Science}, 26:694--706, 1980.

\bibitem{Chankong1983}
V.~Chankong and Y.~Y. Haimes.
\newblock {\em Multiobjective Decision Making: Theory and Methodology}.
\newblock Amsterdam, North-Holland, 1983.

\bibitem{Chuong2018}
T.~D. Chuong.
\newblock Linear matrix inequality conditions and duality for a class of robust
  multiobjective convex polynomial programs.
\newblock {\em SIAM Journal on Optimization}, 28(3):2466--2488, 2018.

\bibitem{Chuong2022}
T.~D. Chuong.
\newblock Second-order cone programming relaxations for a class of
  multiobjective convex polynomial problems.
\newblock {\em Annals of Operations Research}, 311(2):1017--1033, 2022.

\bibitem{Ehrgott2005}
M.~Ehrgott.
\newblock {\em Multicriteria Optimization (2nd ed.)}.
\newblock Springer, Berlin, 2005.

\bibitem{gb08}
M.~Grant and S.~Boyd.
\newblock Graph implementations for nonsmooth convex programs.
\newblock In V.~Blondel, S.~Boyd, and H.~Kimura, editors, {\em Recent Advances
  in Learning and Control}, Lecture Notes in Control and Information Sciences,
  pages 95--110. Springer-Verlag Limited, 2008.
\newblock \url{http://stanford.edu/~boyd/graph_dcp.html}.

\bibitem{cvx}
M.~Grant and S.~Boyd.
\newblock {CVX}: Matlab software for disciplined convex programming, version
  2.1.
\newblock \url{http://cvxr.com/cvx}, Mar. 2014.

\bibitem{GuoJiao2022arXiv}
F.~Guo and L.~G. Jiao.
\newblock A new scheme for approximating the weakly efficient solution set of
  vector rational optimization problems.
\newblock {\em arXiv preprint arXiv.2205.12863}, 2022.

\bibitem{henrion2009gloptipoly}
D.~Henrion, J.-B. Lasserre, and J.~L{\"o}fberg.
\newblock Gloptipoly 3: moments, optimization and semidefinite programming.
\newblock {\em Optimization Methods \& Software}, 24(4-5):761--779, 2009.

\bibitem{Jahn2011}
J.~Jahn.
\newblock {\em Vector Optimization, Theory, Applications, and Extensions
  (2nd)}.
\newblock Springer, Berlin, 2011.

\bibitem{jeyakumar2016convergent}
V.~Jeyakumar, J.~B. Lasserre, G.~Li, and T.-S. Pham.
\newblock Convergent semidefinite programming relaxations for global bilevel
  polynomial optimization problems.
\newblock {\em SIAM Journal on Optimization}, 26(1):753--780, 2016.

\bibitem{Jiao2019}
L.~G. Jiao and J.~H. Lee.
\newblock Finding efficient solutions in robust multiple objective optimization
  with {SOS}-convex polynomial data.
\newblock {\em Annals of Operations Research}, 296(1-2):803--820, 2021.

\bibitem{Jiao2020}
L.~G. Jiao, J.~H. Lee, and Y.~Y. Zhou.
\newblock A hybrid approach for finding efficient solutions in vector
  optimization with {SOS}-convex polynomials.
\newblock {\em Operations Research Letters}, 48(2):188--194, 2020.

\bibitem{kuroiwa2012robust}
D.~Kuroiwa and G.~M. Lee.
\newblock On robust multiobjective optimization.
\newblock {\em Vietnam J. Math}, 40(2-3):305--317, 2012.

\bibitem{Lasserre2010}
J.~B. Lasserre.
\newblock {\em Moments, Positive Polynomials and their Applications}.
\newblock Imperial College Press, London, 2010.

\bibitem{lasserre2010joint+}
J.~B. Lasserre.
\newblock A “joint+ marginal” approach to parametric polynomial
  optimization.
\newblock {\em SIAM Journal on Optimization}, 20(4):1995--2022, 2010.

\bibitem{lasserre2011min}
J.~B. Lasserre.
\newblock Min-max and robust polynomial optimization.
\newblock {\em Journal of Global Optimization}, 51(1):1--10, 2011.

\bibitem{lasserre2012algorithm}
J.-B. Lasserre.
\newblock An algorithm for semi-infinite polynomial optimization.
\newblock {\em Top}, 20(1):119--129, 2012.

\bibitem{Lasserre2015}
J.~B. Lasserre.
\newblock {\em An Introduction to Polynomial and Semi-Algebraic Optimization}.
\newblock Cambridge University Press, 2015.

\bibitem{Lee2018}
J.~H. Lee and L.~G. Jiao.
\newblock Solving fractional multicriteria optimization problems with sum of
  squares convex polynomial data.
\newblock {\em Journal of Optimization Theory and Applications},
  176(2):428--455, 2018.

\bibitem{Lee2021}
J.~H. Lee and L.~G. Jiao.
\newblock Robust multi-objective optimization with {SOS}-convex polynomials
  over a polynomial matrix inequality.
\newblock {\em Journal of Nonlinear and Convex Analysis}, 22(7):1263--1284,
  2021.

\bibitem{Lee2021JOGO}
J.~H. Lee, N.~Sisarat, and L.~G. Jiao.
\newblock Multi-objective convex polynomial optimization and semidefinite
  programming relaxations.
\newblock {\em Journal of Global Optimization}, 80(1):117--138, 2021.

\bibitem{Luc1989}
D.~T. Luc.
\newblock {\em Theory of Vector Optimization}.
\newblock Springer, Berlin, 1989.

\bibitem{nie2013certifying}
J.~Nie.
\newblock Certifying convergence of lasserre’s hierarchy via flat truncation.
\newblock {\em Mathematical Programming}, 142(1):485--510, 2013.

\bibitem{nie2014optimality}
J.~Nie.
\newblock Optimality conditions and finite convergence of lasserre’s
  hierarchy.
\newblock {\em Mathematical programming}, 146(1):97--121, 2014.

\bibitem{nie2015linear}
J.~Nie.
\newblock Linear optimization with cones of moments and nonnegative
  polynomials.
\newblock {\em Mathematical Programming}, 153(1):247--274, 2015.

\bibitem{Putinar1993}
M.~Putinar.
\newblock Positive polynomials on compact semi-algebraic sets.
\newblock {\em Indiana University Mathematics Journal}, 42(3):969--984, 1993.

\bibitem{Sawaragi1985}
Y.~Sawaragi, H.~Nakayama, and T.~Tanino.
\newblock {\em Theory of Multiobjective Optimization}.
\newblock Academic Press Inc., 1985.

\bibitem{Soyster1973}
A.~Soyster.
\newblock Convex programming with set-inclusive constraints and applications to
  inexact linear programming.
\newblock {\em Operations Research}, 21:1154--1157, 1973.

\bibitem{doi:10.1080/10556789908805766}
J.~F. Sturm.
\newblock Using sedumi 1.02, a {MATLAB} toolbox for optimization over symmetric
  cones.
\newblock {\em Optimization Methods and Software}, 11(1-4):625--653, 1999.

\end{thebibliography}
\end{document}